\documentclass[11pt]{article}

\usepackage[hyperref]{package_RFb}
\usepackage{numbered_constants}

\usepackage[numbers]{natbib}
\usepackage[top=3cm,bottom=3cm,left=3.5cm,right=3.5cm]{geometry}

\title{$ \Lambda $-Fleming-Viot processes arising in logistic Bienaymé-Galton-Watson processes with a large carrying capacity}
\author{Rapha\"el Forien}

\DeclareDocumentCommand{\X}{}{\mathcal{X}}
\DeclareDocumentCommand{\MF}{}{\mathcal{M}_F(\X)}
\DeclareDocumentCommand{\M}{}{\mathcal{M}_1(\X)}
\DeclareDocumentCommand{\C}{}{C(\X)}
\DeclareDocumentCommand{\Nbar}{}{\overline{N}_K}
\DeclareDocumentCommand{\Nbeta}{}{\overline{N}_{\beta,K}}

\DeclareConstant{Cst}{C}

\begin{document}
	\maketitle
	
	\begin{abstract}
		We consider a continuous-time Bienaymé-Galton-Watson process with logistic competition in a regime of weak competition, or equivalently of a large carrying capacity. Individuals reproduce at random times independently of each other but die at a rate which increases with the population size. When individuals reproduce, they produce a random number of offspring, drawn according to some probability distribution on the natural integers.
		We keep track of the number of descendants of the initial individuals by adding neutral markers to the individuals, which are inherited by one's offspring.
		We then consider several scaling limits of the measure-valued process describing the distribution of neutral markers in the population, as well as the population size, when the competition parameter tends to zero.
		Three regimes emerge, depending on the tail of the offspring distribution.
		When the offspring distribution admits a second moment (actually a $ 2+\delta $ moment for some positive $ \delta $), the fluctuations of the population size around its carrying capacity are small and the neutral types asymptotically follow a Fleming-Viot process.
		When the offspring distribution has a power-law decay with exponent $ \alpha \in (1,2) $, the population size remains most of the time close to its carrying capacity with some (short-lived) fluctuations, and the neutral types evolve in the limit according to a generalised $ \Lambda $-Fleming-Viot process.
		When the exponent $ \alpha $ is equal to 1, the time scale of the fluctuations changes drastically, as well as the order of magnitude of the population size.
		In that case the limiting dynamics of the neutral markers is given by the dual of the Bolthausen-Sznitman coalescent.
		
		~~
		
	\end{abstract}
	
	\section*{Introduction}
	
	In \cite{kingman_genealogy_1982}, Kingman introduced his famous coalescent to describe the limiting genealogy of a sample of individuals from a Wright-Fisher model as the population size tends to infinity.
	Since then, coalescent theory has flourished in theoretical population genetics and has found many applications thanks to the increasing availability of genetic sequences.
	In order to exploit the information that these sequences carry on the past history of a population, one thus needs to understand how this history shapes the genealogy of individuals randomly sampled from this population.
	
\begin{figure}
		\centering
		\includegraphics[width=0.5\linewidth]{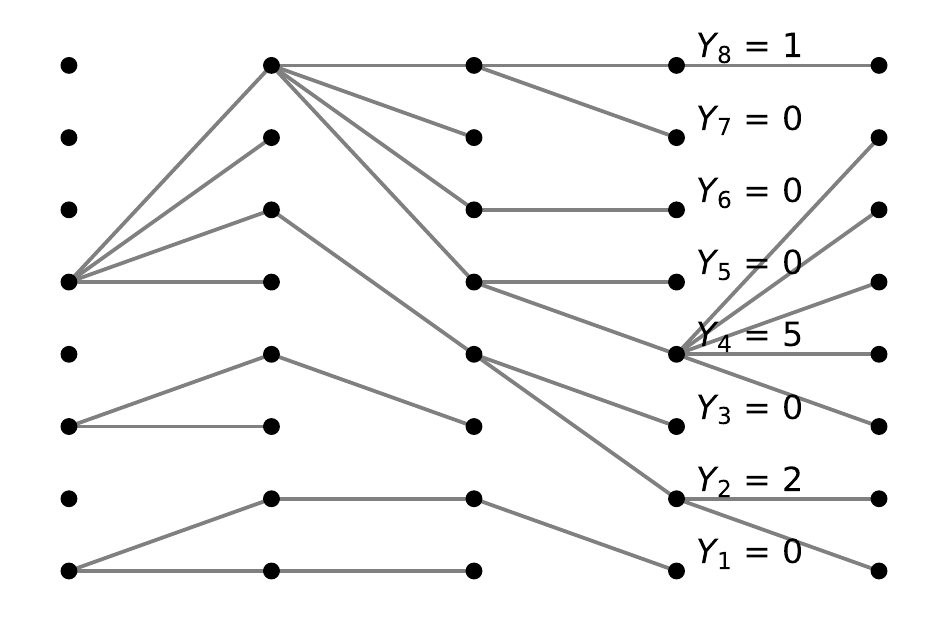}
		\caption{Illustration of a realisation of a Cannings model with $ N = 8 $ individuals, run for 4 generations. The offspring numbers in the last step are indicated on the figure.} \label{fig:Cannings}
	\end{figure}

	Here, we will focus on exchangeable populations, meaning that the individuals composing the population are all equivalent.
	The historic example of such a model is the Cannings model, introduced in \cite{cannings_latent_1974,cannings_latent_1975}, which generalises the original Wright-Fisher model.
	In this model, the population is composed of $ N $ individuals at each (discrete) generation.
	The reproduction mechanism is specified using a random vector $ \bm{Y} = (Y_1, \ldots, Y_N) \in \N^N $ whose distribution is invariant by permutation (i.e. it is exchangeable) and such that
	\begin{equation*}
		\sum_{i=1}^{N} Y_i = N, \quad \text{ almost surely.}
	\end{equation*}
	We then consider a sequence of independent and identically distributed copies of $ \bm{Y} $, denoted by $ (\bm{Y}^{(n)}, n \in \N) $, and the model is defined by stating that the $ i $-th individual in generation $ n $ leaves exactly $ Y_{i}^{(n)} $ descendants in generation $ n+1 $ (see Figure~\ref{fig:Cannings}).
	The classical Wright-Fisher model is obtained by assuming that $ \bm{Y} $ follows a multinomial distribution with uniform weights.
	
	To study the genealogy of the population in this model, we can consider an iid sequence $ (\bm{Y}^{(n)}, n \in \Z) $ and also run it backwards in time.
	We then randomly assign labels $ \lbrace 1, \ldots, N \rbrace $ to the individuals in generation 0, and, for $ n \in \N $, let $ \Pi^N_n $ be the partition of $ \lbrace 1, \ldots, N \rbrace $ induced by the equivalence relation
	\begin{equation*}
		i \sim j \Leftrightarrow \text{ individuals $ i $ and $ j $ in generation 0 share a common ancestor in generation $ -n $.}
	\end{equation*}
	This defines a Markov chain $ (\Pi^N_n, n \in \N) $ taking values in the space of partitions of $ \lbrace 1, \ldots, N \rbrace $.
	For $ 1 \leq m \leq N $, also let $ \Pi^{m,N}_n $ be the partition induced on $ \lbrace 1, \ldots, m \rbrace $ by $ \Pi^N_n $.
	Then $ (\Pi^{m,N}_n, n \in \N) $ is again a Markov chain, and is called a coalescent process (since all its transitions involve merging two or more blocks together).
	
	It can readily be seen that, in the Wright-Fisher model, the probability that any two blocks (also called lineages) merge in the next step is $ 1/N $.
	In fact, in \cite{kingman_genealogy_1982}, Kingman showed that, for this model, for any fixed $ m \in \N $,
	\begin{equation*}
		(\Pi^{m,N}_{\lfloor N t \rfloor}, t \in [0,T]) \cvgas{N} (\Pi^m_t, t \in [0,T]),
	\end{equation*}
	in distribution, where $ (\Pi^m_t, t \geq 0) $ is the continuous-time coalescent process started from the trivial partition $ \lbrace \lbrace 1 \rbrace, \ldots, \lbrace m \rbrace \rbrace $ with transition rates
	\begin{equation*}
		\lambda_{\pi, \pi'} = \begin{cases}
		1 & \text{ if $ \pi' $ is obtained by merging exactly two blocks of $ \pi $,} \\
		0 & \text{ otherwise.}
		\end{cases}
	\end{equation*}
	This process is called the Kingman coalescent (or $ m $-coalescent if one wishes to stress the initial condition).
	
	In \cite{sagitov_general_1999,pitman_coalescents_1999}, Pitman and Sagitov independently introduced a more general family of continuous-time coalescent processes, called $ \Lambda $-coalescents, which allow more than two lineages to merge at the same time, but with no simultaneous mergers.
	The transition rates of a $ \Lambda $-coalescent are characterised by a finite measure on $ [0,1] $, denoted by $ \Lambda $, such that, if $ \pi $ contains $ n $ blocks and $ \pi' $ is obtained by merging exactly $ k $ of them together, then
	\begin{equation*}
		\lambda_{\pi, \pi'} = \int_{0}^{1} u^{k-2} (1-u)^{n-k} \Lambda(du),
	\end{equation*}
	and $ \lambda_{\pi, \pi'} = 0 $ otherwise.
	Note that, when $ \Lambda(du) = \delta_0(du) $, this definition coincides with that of the Kingman coalescent.
	Schweinsberg \cite{schweinsberg_coalescents_2000} and M\"ohle and Sagitov \cite{mohle_classification_2001} then introduced a larger class of coalescents, $ \Xi $-coalescents, which also include simultaneous multiple mergers, and characterised by a finite measure $ \Xi $ on the infinite simplex.
	In addition, \cite{mohle_classification_2001} showed that this family describes all the possible scaling limits of the genealogies of Cannings models as $ N \to \infty $, and gave explicit criteria for such convergence.
	
	Multiple mergers (coalescence events involving more that two lineages) typically arise in population models with skewed offspring distributions.
	In \cite{schweinsberg_coalescent_2003}, Schweinsberg considered a Cannings model in which $ \bm{Y} $ is constructed as follows.
	Consider an $ \N $-valued random variable $ X $ such that $ \E{X} > 1 $, and let $ \bm{X} = (X_1, \ldots, X_N) $ be a family of iid copies of $ X $.
	To obtain $ \bm{Y} $, we then sample $ N $ new individuals without replacement from $ \bm{X} $.
	Equivalently, each individual in the previous generation first produces a random number of propagules distributed as $ X $, independently of the others, and $ N $ randomly chosen propagules mature to form the next generation.
	If $ \sum_{i=1}^{N} X_i < N $, then we may decide on an arbitrary sampling procedure to recover $ N $ individuals in the next generation (this does not affect the result as the probability of this event vanishes as $ N \to \infty $).
	Then let $ (\Pi^{m,N}_n, n \in \N) $ be the coalescent process defined above.
	Schweinsberg then proves the following.
	
	\begin{theorem}[Theorem~4 in \cite{schweinsberg_coalescent_2003}] \label{thm:schweinsberg}
		\begin{enumerate}[i)]
			\item If $ \E{X^2} < \infty $, then $ (\Pi^{m,N}(\lfloor Nt \rfloor), t \geq 0) $ converges in distribution as $ N \to \infty $ to Kingman's coalescent started from $ m $ lineages in which the rate of each pairwise merger is given by
			\begin{equation*}
				\frac{\E{X(X-1)}}{\E{X}^2}.
			\end{equation*}
			\item If $ \P{X \geq k} \sim C k^{-\alpha} $ as $ k \to \infty $ for some constant $ C > 0 $ and $ \alpha \in (1,2) $, then $ (\Pi^{m,N}(\lfloor N^{\alpha-1} t \rfloor), t \geq 0) $ converges in distribution as $ N \to \infty $ to a $ \Lambda $-coalescent started from $ m $ lineages, where $ \Lambda $ is proportional to a $ Beta(2-\alpha, \alpha) $ distribution, i.e.
			\begin{equation*}
				\Lambda(du) = C \alpha \E{X}^{-\alpha} u^{2-\alpha-1}(1-u)^{\alpha-1} du.
			\end{equation*}
			\item If $ \P{X \geq k} \sim C k^{-1} $ as $ k \to \infty $ for some constant $ C > 0 $, then $ (\Pi^{m,N}_{\lfloor \log(N) t \rfloor}, t \in [0,T]) $ converges in distribution as $ N \to \infty $ to a $ \Lambda $-coalescent started from $ m $ lineages, where $ \Lambda $ is Lebesgue measure, i.e.
			\begin{equation*}
				\Lambda(du) = du.
			\end{equation*}
		\end{enumerate}
	\end{theorem}
	
	Other regimes are also treated in \cite{schweinsberg_coalescent_2003}, but they will not be relevant here.
	Note that, in the third part of the statement the limiting coalescent is the Bolthausen-Sznitman coalescent, introduced in \cite{bolthausen_ruelles_1998}.
	
	An informal explanation of Theorem~\ref{thm:schweinsberg} can be given as follows.
	When $ \E{X^2} < \infty $, it is very rare to find an individual that produces $ \bigO{N} $ offspring in a single generation, and the probability that such an event takes place in less than $ \bigO{N} $ generations tends to zero as $ N \to \infty $.
	As a result the family sizes in each generation are all $ \littleO{N} $ with high probability, and hence no multiple mergers take place on this time scale.
	The model thus behaves essentially like a Wright-Fisher model.
	In the second regime, when $ \alpha \in (1,2) $, it take $ \bigO{N^{\alpha-1}} $ generations for an individual to have $ \bigO{N} $ offspring in a single generation, at which point a non-zero fraction of extant lineages can coalesce with some probability.
	The probability that more than one individual have such exceptionally large offspring in a single generation on such time scales is however vanishingly small as $ N \to \infty $, so there are no multiple mergers.
	When $ \alpha = 1 $, the reasoning is similar, except that the typical time until the first multiple merger is of the order of $ \log(N) $.
	
	The Kingman coalescent is also known to satisfy a moment duality with the Wright-Fisher diffusion and more generally with the measure-valued Fleming-Viot diffusion \cite{ethier_flemingviot_1993}, which can be obtained as the limit in distribution as $ N \to \infty $ of the distribution of a (possibly infinite) set of neutral types in the Wright-Fisher model, after scaling time by a factor $ N $.
	The Fleming-Viot diffusion can be characterised as the unique solution to a martingale problem as follows.
	Let $ \X $ be a compact metric space, and, for $ \phi : \X \to \R $ continuous and $ f : \R \to \R $ twice continuously differentiable, let $ F_{f,\phi} : \M \to \R $ be defined on the set $ \M $ of probability measures on $ \X $ (equipped with the weak topology) by
	\begin{equation*}
		F_{f,\phi}(\rho) := f(\langle \rho, \phi \rangle).
	\end{equation*}
	We then define an operator $ \mathcal{L} $ acting on functions of this form as
	\begin{equation*}
		\mathcal{L} F_{f,\phi}(\rho) := \frac{1}{2} f''(\langle \rho, \phi \rangle) \left( \langle \rho, \phi^2 \rangle - \langle \rho, \phi \rangle^2 \right).
	\end{equation*}
	Then, a $ \M $-valued Markov process $ (\rho_t, t \geq 0) $ is a Fleming-Viot diffusion if it satisfies the associated martingale problem, i.e. if for any $ f $ and $ \phi $ as above,
	\begin{equation*}
		F_{f,\phi}(\rho_t) - \int_{0}^{t} \mathcal{L} F_{f,\phi}(\rho_s) ds
	\end{equation*}
	is a martingale with respect to its natural filtration.
	The moment duality with the Kingman coalescent can then be stated as follows (see \cite{ethier_flemingviot_1993}).
	For $ n \geq 1 $, let $ e_n $ denote the partition of $ \lbrace 1, \ldots, n \rbrace $ into singletons, and if $ \pi $ is a partition of $ \lbrace 1, \ldots, n \rbrace $ with $ k $ distinct blocks of respective sizes $ s_1, \ldots, s_k $, we define $ \tau_\pi : \X^k \to \X^n $ as
	\begin{equation*}
		\tau_\pi(y_1, \ldots, y_{k}) = (\underbrace{y_1, \ldots, y_1}_{s_1}, \underbrace{y_2, \ldots, y_2}_{s_2}, \ldots, \underbrace{y_{k}, \ldots, y_{k}}_{s_k}).
	\end{equation*}
	Then, for any continuous function $ \Phi : \X^n \to \R $ invariant by permutation and any $ \rho_0 \in \M $,
	\begin{equation} \label{moment_duality}
		\E[\rho_0]{ \langle \rho_t^{\otimes n}, \Phi \rangle } = \E[e_n]{ \langle \rho_0^{\otimes | \pi_t |}, \Phi \circ \tau_{\pi_t} \rangle },
	\end{equation}
	where $ (\pi_t, t \geq 0) $ is the Kingman coalescent, and $ | \pi | $ denotes the number of blocks of the partition $ \pi $.
	
	Note that this duality implies the well-known duality between the Kingman coalescent and the Wright-Fisher diffusion, since, when $ \X = \lbrace 0, 1 \rbrace $, then
	\begin{equation*}
		\rho_t = X_t \delta_1 + (1-X_t) \delta_0,
	\end{equation*}
	and $ (X_t, t \geq 0) $ is a Wright-Fisher diffusion.
	Choosing $ \Phi(y_1, \ldots, y_n) = \1{y_1 = y_2 = \ldots = y_n = 1} $ in \eqref{moment_duality} yields the classical moment duality
	\begin{equation*}
		\E[X_0]{X_t^n} = \E[e_n]{X_0^{| \pi_t |}}.
	\end{equation*}
	
	This moment duality can be generalised for $ \Lambda $-coalescents, and the relevant measure-valued processes were introduced by Bertoin and Le Gall in \cite{bertoin_stochastic_2003}.
	Given a finite measure $ \Lambda $ on $ [0,1] $, $ \phi : \X \to \R $ continuous and $ f : \R \to \R $ twice continuously differentiable, define an operator $ \mathcal{L}_\Lambda $ acting on functions of the form $ F_{f,\phi} $ as
	\begin{equation} \label{generator_LambdaFV}
		\mathcal{L}_\Lambda F_{f,\phi}(\rho) := \int_{[0,1]} \int_\X \left( f((1-u) \langle \rho, \phi \rangle + u \phi(x)) - f(\langle \rho, \phi \rangle) \right) \rho(dx) \frac{\Lambda(du)}{u^2}.
	\end{equation}
	Bertoin and Le Gall then show that the martingale problem associated to this operator is well posed.
	The corresponding $ \M $-valued Markov process is called a generalised Fleming-Viot process, or a $ \Lambda $-Fleming-Viot process.
	If $ \Lambda(\lbrace 0 \rbrace) = 0 $, then we can construct a $ \Lambda $-Fleming-Viot process $ (\rho_t, t \geq 0) $ with the help of a Poisson random measure $ Q $ on $ \R_+ \times (0,1] $ with intensity $ dt \otimes \frac{\Lambda(du)}{u^2} $, and saying that, for each atom $ (t,u) $ of the random measure $ Q $, a jump takes place in $ \rho $ at time $ t $, in which we sample a type $ x \in \X $ according to $ \rho_{t^-} $, and replace a fraction $ u $ of the population with individuals of type $ x $, so that $ \rho_t $ is obtained by
	\begin{equation*}
		\rho_t = (1-u) \rho_{t^-} + u \delta_x.
	\end{equation*}
	Note that, when $ u \downarrow 0 $, by a Taylor expansion of $ f $ around $ \langle \rho, \phi \rangle $, the integrand in \eqref{generator_LambdaFV} becomes
	\begin{equation*}
		\frac{1}{2} f''(\langle \rho, \phi \rangle) \left( \langle \rho, \phi^2 \rangle - \langle \rho, \phi \rangle^2 \right) + \littleO{1}.
	\end{equation*}
	This ensures that $ \mathcal{L}_\Lambda F_{f,\phi} $ is well defined even if $ \Lambda $ has an atom at 0, and also shows that, when $ \Lambda = \delta_0 $, this definition coincides with that of the Fleming-Viot diffusion.
	The moment duality \eqref{moment_duality} then holds without modification if $ (\rho_t, t \geq 0) $ is a $ \Lambda $-Fleming-Viot process and $ (\pi_t, t \geq 0) $ is a $ \Lambda $-coalescent with the same measure $ \Lambda $ (see Lemma~5 in \cite{bertoin_stochastic_2003}).
	
	All the population models considered so far share a common feature: the population size in each generation is a deterministic constant.
	The main reason for this is that the backwards-in-time dynamics of the genealogy of a sample of individuals depends on the trajectory of the population size, so that if the population is allowed to fluctuate randomly, there is no guarantee that the genealogy of a random sample of individuals can be described by a Markov process running backwards in time (and indeed in many instances it cannot).
	This has made the study of genealogies in populations with fluctuating population sizes especially challenging.
	
	Genealogies of populations with fluctuating population size have nonetheless been the subject of many recent works.
	The case of purely branching populations (without interactions between individuals) has been particularly studied, see \cite{lambert_coalescence_2003,harris_coalescent_2020,harris_universality_2024} and references therein.
	For interacting populations as in the model presented below, however, little is known.
	Several approaches have been proposed to define and characterise the genealogies of interacting populations, such as lookdown constructions, also called countable representations \cite{donnelly_particle_1999,kurtz_poisson_2011,etheridge_genealogical_2019} or tree-valued diffusions \cite{greven_convergence_2009,greven_tree-valued_2013}.
	In general, one can represent the genealogy of a (possibly infinite) population as an ultrametric measured space, as in \cite{foutel-rodier_exchangeable_2021}, and either study the dynamics of this object forwards-in-time or study its distribution at a fixed (typical) time (for example using spinal decomposition techniques \cite{harris_many--few_2017,foutel-rodier_convergence_2023,bansaye_spine_2024}).
	We conjecture that the above techniques could be applied in our setting, but we will here content ourselves with the study of the forwards-in-time dynamics of the measure-valued process describing the evolution of the frequencies of an arbitrary (and possibly infinite) set of neutral types within the population.
	If, as will be the case below, these processes converge to a well-known measure-valued process with a known moment dual, we argue that this strongly indicates that the corresponding genealogies should converge to this moment dual.
	Establishing this would require additional technical notions and is thus beyond the scope of the present paper, but would constitute a very natural extension.
	
	In this paper, we study a family of continuous time Bienaymé-Galton-Watson processes with logistic competition, in which individuals reproduce at random times independently of each other but die at a rate which depends on the current number of individuals in the population.
	The model is defined in such a way that the process is supercritical when there are few individuals but becomes subcritical when the population size exceeds some threshold (called the carrying capacity).
	When individuals reproduce, they instantly produce a random number of new individuals, drawn according to some probability distribution on $ \N $.
	Such processes were introduced in \cite{lambert_branching_2005}.
	
	We study a regime of weak competition, when the carrying capacity of the population tends to infinity.
	More precisely we let the intensity of competition tend to zero as $ 1/K $ for some scaling parameter $ K $ which tends to infinity, so that the population size is typically proportional to $ K $.
	We do not study the genealogy of the population directly, but we are interested in the limiting behaviour of the distribution of an arbitrary set of neutral types in the population.
	This yields a measure-valued process indexed by the scaling parameter $ K $, which tracks the distribution of these neutral markers in the population.
	
	We study three regimes, determined by the tail of the offspring distribution.
	In each regime, we show that, on a suitable time scale, the population size, suitably rescaled, spends most of its time close to a deterministic constant, and that, on the same time scale, the distribution of neutral markers converges to a $ \Lambda $-Fleming-Viot process defined above.
	
	We consider three possible regimes, corresponding to the three cases of Theorem~\ref{thm:schweinsberg} above.
	When the number of offspring produced at each reproduction admits a moment of order $ 2 + \delta $ for some $ \delta > 0 $, we show that the limiting dynamics of neutral types is given by a classical Fleming-Viot process, and that this convergence takes place on time intervals of the form $ [0, K T] $.
	Moreover, with high probability the population size remains in sets of the form $ [K(n_*-\varepsilon), K(n_* + \varepsilon)] $ on this time interval for some explicit $ n_* > 0 $, for any $ \varepsilon >0 $ (see Theorem~\ref{thm:finite_variance} below).
	If the reproduction law has a power-law tail with exponent $ \alpha \in (1,2) $ as in case ii) of Theorem~\ref{thm:schweinsberg}, then the limiting dynamics of neutral markers is the dual of a $ \mathrm{Beta}(2-\alpha, \alpha) $-coalescent, and this convergence takes place on time intervals of the form $ [0, K^{\alpha-1} T] $.
	Moreover, on such time intervals, the amount of time that the population size spends outside sets of the form $ [K(n_* - \varepsilon), K(n_* + \varepsilon)] $ tends to zero as $ K $ tends to infinity for any $ \varepsilon >0 $ (see Theorem~\ref{thm:alpha-stable} below).
	In the case $ \alpha = 1 $, the dynamics of neutral markers converge to the dual of a Bolthausen-Sznitman coalescent on time intervals of the form $ [0,T] $, and, on such time intervals, the amount of time that the population size spends outside sets of the form $ [K \log(K) (n_0 - \varepsilon), K\log(K) (n_0 + \varepsilon)] $ tends to zero as $ K \to \infty $, for any $ \varepsilon > 0 $, for some explicit $ n_0 > 0 $ (Theorem~\ref{thm:neveu} below).
	
	We note that these time scales are highly reminiscent of the ones appearing in Theorem~\ref{thm:schweinsberg}.
	Even in the last regime ($ \alpha = 1 $), the ratio between the population size and the time scale over which the dynamics of neutral markers converges is the same as in case iii of Theorem~\ref{thm:schweinsberg}.
	We note that our result in the finite variance case is a generalisation of a similar result in \cite{billiard_stochastic_2015} in the case of purely binary branching, see Remark~\ref{rk:FV-billiard} below.
	
	The results obtained in the present work are in direct correspondence with those of \cite{birkner_alpha-stable_2005}, in which the frequency process associated to $ \alpha $-stable continuous-state branching processes was studied.
	In particular, after a time change depending on $ \alpha $ and on the total population size, the frequency process reduces to a $ \Lambda $-Fleming-Viot process.
	Such $ \alpha $-stable continuous-state branching processes can be obtained as scaling limits of the continuous-time Bienaymé-Galton-Watson processes that we study here, in precisely the same regime as our main results.
	The additional effect of competition is to force the (scaled) population size to stay close to a fixed value, without affecting the frequency process.
	This makes the time change that appears in the results of \cite{birkner_alpha-stable_2005} trivial, which yields the three time scales given in the results below.
	
	The rest of the paper is structured as follows.
	In Section~\ref{sec:model_results}, we introduce the model and state the main results.
	Each of the three following sections is then dedicated to the proof of one of the three main theorems.

	\section{Definition of the model and main results} \label{sec:model_results}
	
	We consider a population of individuals where each individual carries a neutral label belonging to some compact metric space $ \X $.
	Each individual, at some rate $ b > 0 $, produces a random number of children chosen according to some probability distribution $ (p_k, k \in \N \setminus \lbrace 0 \rbrace) $.
	Each child inherits the label of their parent.
	In addition, each individual dies at an instantaneous rate $ d + c_K N(t) $, where $ d > 0 $, $ c_K > 0 $ and $ N(t) $ is the number of individuals alive at time $ t \geq 0 $.
	
	The state of the population at time $ t \geq 0 $ is then represented by a finite point measure $ \nu^K_t $ on $ \X $ where, if $ (X_1(t), \ldots, X_{N(t)}(t)) $ is the collection of labels of the individuals alive at time $ t $,
	\begin{equation*}
		\nu^K_t := \sum_{i=1}^{N(t)} \delta_{X_i(t)}.
	\end{equation*}
	Let $ h : [0,1] \to [0,1] $ be the generating function associated to the probability distribution $ (p_k, k \geq 1) $, i.e.
	\begin{equation*}
		h(s) := \sum_{k=1}^{+\infty} s^k p_k, \quad s \in [0,1].
	\end{equation*}
	We then assume that, for any $ \varepsilon \in (0,1) $,
	\begin{equation} \label{assumption_h}
		\int_{1-\varepsilon}^{1} \frac{1}{| h(s) - s |} ds = \infty.
	\end{equation}
	By Theorem~V.9.1 in \cite{harris_theory_1963}, this ensures that the process is well defined for all time $ t \geq 0 $ almost surely, as it is stochastically dominated by a non-explosive pure birth branching process.
	Note that \eqref{assumption_h} is satisfied as soon as $ h'(1) = \sum_{k=1}^{\infty} k p_k $ is finite.

	Let $ \MF $ denote the space of finite measures on $ \X $, equipped with the topology of weak convergence and let $ \C $ denote the space of continuous (hence bounded) real-valued functions on $ \X $.
	For any $ \nu \in \MF $ and $ \phi \in \C $, let $ \langle \nu, \phi \rangle $ denote the integral
	\begin{equation*}
		\langle \nu, \phi \rangle := \int_{\X} \phi \, d\nu.
	\end{equation*}
	For $ f : \R \to \R $ smooth and bounded and $ \phi \in \C $, define $ F_{f,\phi} : \MF \to \R $ as
	\begin{equation*}
		F_{f,\phi}(\nu) := f(\langle \nu, \phi \rangle), \quad \nu \in \MF.
	\end{equation*}
	Let $ \mathcal{D}_0 $ denote the space of real-valued functions on $ \MF $ of this form, and let $ \mathcal{D} $ denote its linear span.
	By the Stone-Weierstrass theorem, $ \mathcal{D} $ is dense in $ C(\M) $ as it contains the set of linear combinations of functions of the form $ \nu \mapsto e^{-\langle \nu, \phi \rangle} $, $ \phi \in C(\X) $, which is a separating subalgebra.
	Then, $ (\nu^K_t, t \geq 0) $ is a Markov process which solves the martingale problem associated to the generator $ \mathcal{L}_K $ defined on $ \mathcal{D} $ as
	\begin{multline} \label{def:generator}
		\mathcal{L}_K F_{f,\phi}(\nu) := b \sum_{k=1}^{+\infty} p_k \int_\X \left( f(\langle \nu, \phi \rangle + k \phi(x)) - f(\langle \nu, \phi \rangle) \right) \nu(dx) \\ + (d + c_K \langle \nu, 1 \rangle) \int_\X \left( f(\langle \nu, \phi \rangle - \phi(x)) - f(\langle \nu, \phi \rangle) \right) \nu(dx).
	\end{multline}
	Note that the labels carried by the individuals are neutral, since they do not affect their reproductive success in any way.
	As a result the process $ t \mapsto \langle \nu^K_t, 1 \rangle $ is also a Markov process, and we can see that it is distributed as a so-called logistic birth and death process \cite{lambert_branching_2005}.
	The neutral labels should be seen as a way to keep track of the fate of the families descended from the individuals alive at time zero (where by family we mean individuals who carry the same label, as they are necessarily descended from individuals who carried the same label at time zero).
	
	The parameter $ K $ is a scaling parameter that is used to let $ c_K $ tend to zero.
	We thus assume in all the following that there exists $ c > 0 $ such that
	\begin{equation*}
		c_K = \frac{c}{K}.
	\end{equation*}
	We are interested in the behaviour of $ (\nu^K_t, t \geq 0) $ as $ K \to \infty $, when started with a large number of individuals at time 0.
	If $ b m < d $, where
	\begin{equation} \label{def:m}
		m := \sum_{k = 1}^{+\infty} k p_k,
	\end{equation}
	then the number of individuals is stochastically dominated by a subcritical branching process, and thus will go extinct in a time $ \bigO{\log(K)} $ if started from $ \bigO{K} $ individuals at time 0.
	We thus assume in all what follows that $ bm > d $ (but not necessarily that $ m < \infty $).
	We will see that, when $ m < \infty $, the number of individuals in the population will remain most of the time in the neighbourhood of $ K n_* $, where
	\begin{equation} \label{def:n_star}
		n_* := \frac{bm - d}{c}.
	\end{equation}
	(We will make this statement more precise in the following, in particular, we should note that the population will reach extinction in finite time almost surely for any value of $ K $, but the time it takes for this to happen is astronomic compared to the timescales in which we are interested here.)
	
	Along with the number of individuals in the population, we will be interested in the evolution of the distribution of the neutral labels in the population.
	To do so, we define, for any $ \nu \in \MF $, a probability measure $ \varrho(\nu) $ on $ \X $ by
	\begin{equation} \label{def:varrho}
		\varrho(\nu) = \begin{cases}
			\frac{\nu}{\langle \nu, 1 \rangle}, & \text{ if } \langle \nu, 1 \rangle > 0, \\
			\delta_{x_0} & \text{ otherwise,}
		\end{cases}
	\end{equation}
	where $ x_0 $ is an arbitrary element of $ \X $ (for example one such that $ \nu^K_0(\lbrace x_0 \rbrace) = 0 $ almost surely, but the choice of $ x_0 $ will be irrelevant in the following).
	Then let $ \M $ denote the space of probability measures on $ \X $, also equipped with the topology of weak convergence.
	
	We will study the behaviour of $ (\nu^K_t, t \geq 0) $ in three different regimes: the finite variance regime, in which we will assume that $ \sum_{k=1}^{+\infty} k^{2+\delta} p_k < +\infty $ for some $ \delta > 0 $, and the $ \alpha $-stable regime, in which 
	\begin{equation*}
		p_k \propto \frac{1}{k^{1+\alpha}} \quad \text{ as } k \to \infty
	\end{equation*}
	for some $ \alpha \in (1,2) $ on the one hand, and for $ \alpha = 1 $ on the other hand.
	
	\subsection{The finite variance regime} \label{subsec:finite_variance}
	
	In this regime, the relevant time scale is $ [0,KT] $, for $ T > 0 $, and the number of individuals in the population will typically be of the order of $ K $.
	We thus set, for $ t \geq 0 $,
	\begin{align*}
		\Nbar(t) := \frac{1}{K} \langle \nu^K_{K t}, 1 \rangle, && \rho^K_t := \varrho(\nu^K_{K t}).
	\end{align*}
	For a (locally) compact metric space $ S $, we denote by $ D([0,T], S) $ the space of \cadlag $ S $-valued processes indexed by $ [0,T] $, equipped with the usual Skorokhod topology and the J1 metric, with respect to which it is complete and separable \citep[Theorem~12.2]{billingsley_convergence_1999}.
	The main result is the following.
	
	\begin{theorem} \label{thm:finite_variance}
		Assume that there exists $ \delta > 0 $ such that
		\begin{equation} \label{assumption_p_finite_variance}
			\sum_{k=1}^{+\infty} k^{2+\delta} p_k < +\infty.
		\end{equation}
		Then set
		\begin{align*}
			m := \sum_{k=1}^{+\infty} k p_k, && m_{(2)} := \sum_{k=1}^{+\infty} k^2 p_k.
		\end{align*}
		Assume that $ bm > d $ and set $ n_* := \frac{bm - d}{c} $.
		Also assume that $ \Nbar(0) \to n_* $ and $ \rho^K_0 \to \rho_0 $ in probability as $ K \to \infty $ for some $ \rho_0 \in \M $.
		Then, for any $ \varepsilon > 0 $ and any $ T > 0 $,
		\begin{equation} \label{FW_bound_pop_size_finite_var}
			\lim_{K \to \infty} \P{ \sup_{t \in [0,T]} | \Nbar(t) - n_* | > \varepsilon } = 0.
		\end{equation}
		In addition, as $ K \to \infty $, $ (\rho^K_t, t \in [0,T]) $ converges in distribution in $ D([0,T], \M) $ to a Fleming-Viot diffusion on $ \X $ with speed $ 1/N_e $, where
		\begin{equation} \label{effective_pop_size}
			N_e := \frac{n_*}{b(m + m_{(2)})}.
		\end{equation}
	\end{theorem}
	
	We note that we assume that $ \Nbar(0) $ converges to the equilibrium value $ n_* $ mostly for convenience, it would be enough to assume that it converges to any other positive value, or even that it is tight in $ (0,\infty) $, and the statement of the above result would hold with only the slight modification that the supremum over $ t $ in \eqref{FW_bound_pop_size_finite_var} should be taken over $ [\delta, T] $ for some arbitrarily small $ \delta > 0 $.
	This is easily seen by using the fact that, if $ \Nbar(0) $ converges to some other value $ n_0 \in (0,\infty) $, then, by the Law of large numbers, $ (\Nbar(t/K), t \in [0,T]) $ remains in a neighbourhood of the trajectory of the solution to the deterministic logistic equation with high probability as $ K \to \infty $, for any $ T > 0 $.
	Since this deterministic solution, started from any positive initial condition, converges to $ n_* $ as $ t \to \infty $, we can find $ T $ such that $ \Nbar(T/K) $ is in a small neighbourhood of $ n_* $ with high probability as $ K \to \infty $.
	
	\begin{remark} \label{rk:FV-billiard}
		Theorem~\ref{thm:finite_variance} generalises a result of \cite{billiard_stochastic_2015} (Theorem~3.4 and Proposition~4.7), in which a similar model was considered with $ p_1 = 1 $ (and $ p_k = 0 $ for all $ k \geq 2 $).
		The authors consider a model for a population of individuals carrying both a trait under selection and a neutral marker, under the assumption that mutations affecting the trait under selection take place much less frequently than those affecting the neutral trait.
		They show that, on a suitable time scale, the dynamics of the distribution of traits and neutral markers is described by a trait substitution sequence for the trait under selection, superimposed with a Fleming-Viot diffusion for the neutral marker, whose speed depends on the current resident trait.
		In our notations, their speed is given by $ \frac{n_*}{2 b} $, where $ n_* $ and $ b $ are the (scaled) population size and birth rate of the current resident population.
		This coincides with \eqref{effective_pop_size}, since $ p_1 = 1 $ yields $ m = m_{(2)} = 1 $.
	\end{remark}

	\subsection[The alpha-stable regime]{The $ \alpha $-stable regime} \label{subsec:alpha-stable}
	
	In this regime, the relevant time scale is $ [0,K^{\alpha-1} T] $, and the typical population size will again be of the order of $ K $.
	We thus set, for $ t \geq 0 $,
	\begin{align*}
		\Nbar(t) := \frac{1}{K} \langle \nu^K_{K^{\alpha-1} t}, 1 \rangle, && \rho^K_t := \varrho(\nu^K_{K^{\alpha-1} t}).
	\end{align*}
	We then have the following result.
	
	\begin{theorem} \label{thm:alpha-stable}
		Assume that there exist $ p_0 > 0 $ and $ \alpha \in (1,2) $ such that
		\begin{equation} \label{assumption:p_k}
			p_k \sim \frac{p_0}{k^{1+\alpha}}, \quad \text{ as } k \to \infty.
		\end{equation}
		Also assume that $ (\Nbar(0), K \geq 1) $ is tight in $ (0,\infty) $ and that $ \rho^K_0 \to \rho_0 $ in probability as $ K \to \infty $ for some $ \rho_0 \in \M $. 
		Then, for any $ \varepsilon > 0 $ and $ T > 0 $,
		\begin{equation} \label{bound_pop_size_alpha-stable}
			\lim_{K \to \infty} \P{ \int_{0}^{T} \1{| \Nbar(t) - n_* | > \varepsilon} dt > \varepsilon} = 0.
		\end{equation}
		In addition, as $ K \to \infty $, $ (\rho^K_t, t \in [0,T]) $ converges in distribution in $ D([0,T], \M) $ to a $ \Lambda $-Fleming-Viot process where
		\begin{equation} \label{Lambda_measure}
			\Lambda(du) = \frac{b\, p_0\, \alpha}{n_*^{\alpha-1}} (1-u)^{\alpha - 1} u^{1-\alpha} du.
		\end{equation}
	\end{theorem}
	
	Theorem~\ref{thm:alpha-stable} strongly parallels point \textit{ii} of Theorem~\ref{thm:schweinsberg}, since the limit we obtain for the distribution of neutral types is the dual of the limiting genealogy obtained there.
	This can be explained by the fact that, since the population size is roughly constant (and of the order of $ K $), the rate at which a single individual produces at least $ y K $ children is roughly 
	\begin{equation*}
		b n_* K \sum_{k \geq y K} p_k \sim \frac{b n_* p_0}{\alpha} K^{1-\alpha} y^{-\alpha}.
	\end{equation*}
	In Schweinsberg's model, the probability that a single individual produces at least $ y N $ children before the sampling step is also proportional to $ N^{1-\alpha} y^{-\alpha} $.
	In our model, competition quickly brings back the population size to a neighbourhood of $ n_* K $, and does so by affecting all types equally, so that this phase effectively acts as the sampling step in Schweinsberg's model.
	
	Theorem~\ref{thm:alpha-stable} is proved in Section~\ref{sec:alpha_stable}, using the stochastic averaging principle of T. Kurtz \cite{kurtz_averaging_1992}.
	Here, the fast process is $ (\Nbar(t), t \geq 0) $ and the slow process is $ (\rho^K_t, t \geq 0) $.

	\subsection{Neveu's branching process with logistic competition} \label{subsec:neveu}
	
	When $ \alpha = 1 $ in the previous regime, the behaviour of $ \Nbar $ becomes degenerate as $ m = +\infty $, so $ n_* = +\infty $ as well.
	On the other hand, it is known that, for the corresponding continuous state branching process without competition, the distribution of types $ \varrho(\nu^K_t) $ is a Markov process, independent of the total population size, up to the extinction time of the process, and that it evolves according to the dual of the Bolthauzen-Sznitman coalescent, i.e. a $ \Lambda $-Fleming-Viot process with $ \Lambda(du) \propto du $.
	
	In this case, we set, for $ t \geq 0 $
	\begin{align*}
		\Nbar(t) := \frac{1}{K \log(K)} \langle \nu^K_{t}, 1 \rangle, && \rho^K_t := \varrho(\nu^K_t).
	\end{align*}
	We then have the following result.
	
	\begin{theorem} \label{thm:neveu}
		Assume that there exists $ p_0 > 0 $ such that
		\begin{equation} \label{condition_pk_neveu}
			\lim_{k \to \infty} k^2 p_k = p_0,
		\end{equation}
		and set
		\begin{equation} \label{def:nbar}
			n_* := \frac{b p_0}{c}.
		\end{equation}
		Also assume that $ (\Nbar(0), K \geq 1) $ is tight in $ (0,\infty) $ and that $ \rho^K_0 \to \rho_0 $ in probability as $ K \to \infty $ for some $ \rho_0 \in \M $. 
		Then, for any $ \varepsilon > 0 $ and $ T > 0 $,
		\begin{equation} \label{concentration_Nbar_neveu}
			\lim_{K \to \infty} \P{ \int_{0}^{T} \1{\abs{\Nbar(s) - n_*} > \varepsilon} ds > \varepsilon } = 0.
		\end{equation}
		In addition, as $ K \to \infty $, $ (\rho^K_t, t \in [0,T]) $ converges in distribution in $ D([0,T], \M) $ to a $ \Lambda $-Fleming-Viot process with
		\begin{equation*}
			\Lambda(du) = b p_0 du.
		\end{equation*}
	\end{theorem}
	
	Note that \eqref{condition_pk_neveu} implies the non explosion criterion \eqref{assumption_h}.
	To see this, we note that, setting $ s = 1-u $,
	\begin{align*}
		s - h(s) &= 1-u - h(1-u) \\
		&= \sum_{k=1}^{+\infty} p_k (1-u) f_k(u), 
	\end{align*}
	where $ f_k(u) = 1-(1-u)^k $.
	Since $ f_k $ is concave,
	\begin{equation*}
		f_k(u) \leq f_k'(0) u \leq (k-1) u.
	\end{equation*}
	Using this when $ k \leq 1/u $, and $ f_k \leq 1 $ otherwise, we obtain
	\begin{equation*}
		1-u - h(1-u) \leq u(1-u) \sum_{k=1}^{\lfloor 1/u \rfloor} k p_k + (1-u) \sum_{k > 1/u} p_k.
	\end{equation*}
	We then note that there exists $ C > 0 $ such that, for $ u $ small enough,
	\begin{align*}
		\sum_{k=1}^{\lfloor 1/u \rfloor} k p_k \leq C \log(1/u), && \text{ and } && \sum_{k > 1/u} p_k \leq C u.
	\end{align*}
	It follows that, for $ u $ small enough,
	\begin{equation*}
		0 \leq 1-u - h(1-u) \leq C(1 + \log(1/u)) u(1-u),
	\end{equation*}
	which implies \eqref{assumption_h}.
	
	At first sight, Theorem~\ref{thm:neveu} might appear to differ significantly from \textit{iii} in Theorem~\ref{thm:schweinsberg}.
	The main reason is that, here the population size is not as constrained as in \cite{schweinsberg_coalescent_2003}, and the $ K \log(K) $ scale appears when we look for a balance between competition and the frequent births (ignoring the birth events producing $ \bigO{K} $ children, as these will happen less frequently).
	To see that the above regime corresponds to \textit{iii} in Theorem~\ref{thm:schweinsberg}, we can note that the \emph{product} of the typical number of individuals and the time scale of the fluctuations of the genetic types is the same in both results.
	
	Theorem~\ref{thm:neveu} is proved in Section~\ref{sec:neveu}, also using Kurtz' averaging principle.
	The tightness condition on the population size process $ (\Nbar(t), t \geq 0) $ is more delicate in this case as the fluctuations do not have a finite first moment, but the condition can be checked by considering a sufficiently slowly growing function, see Subsection~\ref{subsec:Nbar_neveu}.
	
	\section{Convergence to a Fleming-Viot diffusion in the finite variance case} \label{sec:proof_finite_variance}
	
	\subsection{Proof of the main result} \label{subsec:proof_thm_finite_variance}
	
	In this section, we prove Theorem~\ref{thm:finite_variance}.
	Recall that in this case, for $ t \geq 0 $,
	\begin{align*}
		\Nbar(t) = \frac{1}{K} \langle \nu^K_{K t}, 1 \rangle, && \rho^K_t = \varrho(\nu^K_{K t}).
	\end{align*}
	Also let $ (\F^K_t, t \geq 0) $ denote the natural filtration associated to $ (\nu^K_{Kt}, t \geq 0) $ (or equivalently to $ ((\Nbar(t), \rho^K_t), t \geq 0) $).
	We then have the following result, which already yields the first part of the statement of Theorem~\ref{thm:finite_variance}, i.e. \eqref{FW_bound_pop_size_finite_var}.
	
	\begin{lemma} \label{lemma:FW_finite_variance}
		Under the assumptions of Theorem~\ref{thm:finite_variance}, for any $ \varepsilon > 0 $,
		\begin{equation*}
			\lim_{K \to \infty} \P{ \sup_{t \in [0,T]} \abs{ \Nbar(t) - n_* } > \varepsilon } = 0.
		\end{equation*}
	\end{lemma}

	We prove this lemma in Subsection~\ref{subsec:exit_time}.
	In the rest of this subsection, we prove the second part of the statement of Theorem~\ref{thm:finite_variance}, that is the convergence of $ (\rho^K_t, t \in [0,T]) $ to a Fleming-Viot diffusion as $ K \to \infty $.
	The following lemma is key to proving this convergence.
	For $ k \in \N $, let $ \mathcal{C}^k(\X) $ denote the space of real-valued functions on $ \X $ which admit continuous derivatives up to order $ k $.
	
	\begin{lemma} \label{lemma:rescaled_generator_rho}
		Under the assumptions of Theorem~\ref{thm:finite_variance}, for any $ F_{f,\phi} \in \mathcal{D}_0 $ with $ f \in \mathcal{C}^3(\X) $ and for any $ \nu \in \MF $ with $ \langle \nu, 1 \rangle > 1 $ and $ \varrho(\nu) = \rho $,
		\begin{equation} \label{generator_rho}
			K \mathcal{L}_K (F_{f,\phi} \circ \varrho)(\nu) = \frac{1}{2} \mathcal{S}_K(\langle \nu, 1 \rangle) f''(\langle \rho, \phi \rangle)  \int_\X \left( \phi(x) - \langle \rho, \phi \rangle \right)^2 \rho(dx) + G^K_{f,\phi}(\nu),
		\end{equation}
		where, for $ N > 1 $,
		\begin{equation} \label{def:SK}
			\mathcal{S}_K(N) := b \sum_{k=1}^{+\infty} N K \left( \frac{k}{N + k} \right)^2 p_k + \left( d + \frac{c}{K} N \right) \frac{N K}{(N-1)^2},
		\end{equation}
		and where, setting $ N = \langle \nu, 1 \rangle $,
		\begin{equation} \label{bound_GK}
			\abs{ G^K_{f,\phi}(\nu) } \leq \frac{4}{3} \| f^{(3)} \|_{\infty} \| \phi \|^3 \left( b \sum_{k=1}^{+\infty} N K \left( \frac{k}{N+k} \right)^3 p_k + \left( d + \frac{c}{K} N \right) \frac{N K}{(N-1)^3} \right).
		\end{equation}
	\end{lemma}
	
	To see why Theorem~\ref{thm:finite_variance} follows from this, note that
	\begin{align*}
		\mathcal{S}_K(K n_*) &= b n_* \sum_{k=1}^{+\infty} \frac{k^2}{(n_* + \frac{k}{K})^2} p_k + (d + c n_*) n_* \frac{1}{(n_* - \frac{1}{K})^2} \\
		&\cvgas{K} \frac{b}{n_*} m_{(2)} + (d + c n_*) \frac{1}{n_*} = \frac{1}{N_e}. \numberthis \label{cvg_SK}
	\end{align*}
	Moreover, we can use \eqref{bound_GK} to show that $ \abs{ G^K_{f,\phi}(\nu) } = \littleO{1} $ when $ \langle \nu, 1 \rangle = \bigO{K} $.
	As a result, we can reasonably expect that the limit of any converging subsequence of $ (\rho^K_t, t\in [0,T]) $ solves the martingale problem associated to the generator $ (\mathcal{L}, \mathcal{D}) $, where, for $ F = F_{f,\phi} \in \mathcal{D}_0 $,
	\begin{equation*}
		\mathcal{L} F_{f,\phi}(\rho) = \frac{1}{2 N_e} f''(\langle \rho, \phi \rangle) \int_\X \left( \phi(x) - \langle \rho, \phi \rangle \right)^2 \rho(dx),
	\end{equation*}
	which is the generator of the Fleming-Viot diffusion on $ \X $ with speed $ 1/N_e $.
	
	\begin{proof}[Proof of Lemma~\ref{lemma:rescaled_generator_rho}]
		By \eqref{def:generator}, and the definition of $ \varrho(\nu) $ in \eqref{def:varrho},
		\begin{multline} \label{generator_F_rho}
			K \mathcal{L}_K (F_{f,\phi} \circ \varrho) (\nu) = K b \langle \nu, 1 \rangle \sum_{k=1}^{+\infty} p_k \int_\X \left( f(\langle \varrho(\nu + k \delta_x), \phi \rangle) - f(\langle \varrho(\nu), \phi \rangle) \right) \varrho(\nu)(dx) \\ + K (d + c_K \langle \nu, 1 \rangle) \langle \nu, 1 \rangle \int_\X \left( f(\langle \varrho(\nu - \delta_x), \phi \rangle) - f(\langle \varrho(\nu), \phi \rangle) \right) \varrho(\nu)(dx).
		\end{multline}
		We then compute, for $ k \in \N \cup \lbrace -1 \rbrace $,
		\begin{align*}
			\langle \varrho(\nu + k \delta_x), \phi \rangle &= \frac{N \langle \rho, \phi \rangle + k \phi(x)}{N + k} \\
			&= \langle \rho, \phi \rangle + \frac{k}{N + k} \left( \phi(x) - \langle \rho, \phi \rangle \right), \numberthis \label{jump_rho}
		\end{align*}
		where $ \rho = \varrho(\nu) $ and $ N = \langle \nu, 1 \rangle > 1 $.
		Then, by Taylor's formula,
		\begin{multline} \label{taylor_f}
			f\left( \frac{N \langle \rho, \phi \rangle + k \phi(x)}{N + k} \right) = f(\langle \rho, \phi \rangle) + f'(\langle \rho, \phi \rangle) \frac{k}{N + k} \left( \phi(x) - \langle \rho, \phi \rangle \right) \\ + \frac{1}{2} f''(\langle \rho, \phi \rangle) \left( \frac{k}{N + k} \right)^2 \left( \phi(x) - \langle \rho, \phi \rangle \right)^2 + R_{k,x}(N, \rho) \left( \frac{k}{N + k} \right)^3 \left( \phi(x) - \langle \rho, \phi \rangle \right)^3,
		\end{multline}
		where
		\begin{equation*}
			R_{k,x}(N, \rho) := \frac{1}{2} \int_{0}^{1} (1-u)^2 f^{(3)}\left( \langle \rho, \phi \rangle + u \frac{k}{N + k} \left( \phi(x) - \langle \rho, \phi \rangle \right) \right) du.
		\end{equation*}
		Plugging this in \eqref{generator_F_rho}, we note that the factor in front of $ f'(\langle \rho, \phi \rangle) $ cancels when integrated against $ \rho(dx) $.
		A straightforward computation then yields \eqref{generator_rho}, with
		\begin{multline*}
			G^K_{f,\phi}(\nu) = K b \sum_{i=1}^{+\infty} p_k \left( \frac{k}{N + k} \right)^3 N \int_{\X} R_{k,x}(N, \rho) \left( \phi(x) - \langle \rho, \phi \rangle \right)^3 \rho(dx) \\ + K \left( d + \frac{c}{K} N \right) N \left( \frac{-1}{N - 1} \right)^3 \int_{\X} R_{-1,x}(N, \rho) \left( \phi(x) - \langle \rho, \phi \rangle \right)^3 \rho(dx).
		\end{multline*}
		Noting that $ | R_{x,k}(N,\rho) | \leq \frac{1}{6} \| f^{(3)} \|_{\infty} $, we obtain \eqref{bound_GK}.
	\end{proof}
	
	Let us now conclude the proof of Theorem~\ref{thm:finite_variance} by making the argument above rigorous.
	Fix $ \varepsilon > 0 $ and let $ \tau_K $ be the stopping time defined by
	\begin{equation} \label{def:tau_K}
		\tau_K := \inf \lbrace t \geq 0 : \abs{\Nbar(t) - n_*} > \varepsilon \rbrace.
	\end{equation}
	By Lemma~\ref{lemma:FW_finite_variance}, we then have $ \P{\tau_K \leq T} \to 0 $ as $ K \to \infty $.
	
	\begin{lemma} \label{lemma:tightness_finite_var}
		Under the assumptions of Theorem~\ref{thm:finite_variance}, the sequence of $ D([0,T], \M) $-valued random variables $ (\rho^K, K \geq 1) $ is tight.
	\end{lemma}
	
	\begin{proof}
		First, by Lemma~\ref{lemma:FW_finite_variance}, it is enough to show that $ (\rho^K_{\cdot \wedge \tau_K}, K \geq 1) $ is tight.
		In addition, by \citep[Theorem~2.1]{roelly-coppoletta_criterion_1986}, it suffices to show that, for each $ \phi \in \C $, the sequence of $ D([0,T], \R) $-valued random variables $ (\langle \rho^K_{\cdot \wedge \tau_K}, \phi \rangle, K \geq 1 ) $ is tight.
		By Lemma~\ref{lemma:rescaled_generator_rho}, for each $ \phi \in \C $,
		\begin{equation*}
			t \mapsto \langle \rho^K_{t \wedge \tau_K}, \phi \rangle
		\end{equation*}
		is a square integrable martingale with predictable variation process
		\begin{equation*}
			t \mapsto \int_{0}^{t \wedge \tau_K} \mathcal{S}_K(K \Nbar(s)) \int_X \left( \phi(x) - \langle \rho^K_s, \phi \rangle \right)^2 \rho^K_s(dx) ds.
		\end{equation*}
		We then note that
		\begin{equation*}
			\int_X \left( \phi(x) - \langle \rho^K_s, \phi \rangle \right)^2 \rho^K_s(dx) \leq 2 \| \phi \|^2,
		\end{equation*}
		and, for all $ s \in [0, T \wedge \tau_K] $, for $ K $ large enough that $ 1/K \leq \varepsilon $,
		\begin{equation*}
			\mathcal{S}_K(K \Nbar(s)) \leq \frac{b m_{(2)} + d + c (n_* + \varepsilon)}{n_* - 2 \varepsilon} (n_* + \varepsilon),
		\end{equation*}
		almost surely.
		Tightness of $ (\langle \rho^K_{\cdot \wedge \tau_K}, \phi \rangle, K \geq 1) $ then follows by a standard convergence criterion, for example Theorem~1 in \cite{aldous_stopping_1978}, which concludes the proof of the lemma.
	\end{proof}
	
	We are now left with showing that the limit of any converging subsequence of $ (\rho^K, K \geq 1) $ solves the martingale problem associated to $ (\mathcal{L}, \mathcal{D}) $.
	By \citep[Theorem~4.8.10]{ethier_markov_1986}, this follows if we prove that, for all $ F \in \mathcal{D} $ and $ h_1, \ldots, h_k $ continuous and bounded functions on $ \M $,
	\begin{equation} \label{cvg_martingale_pb}
		\lim_{K \to \infty} \E{ \left( F(\rho^K_{t+s}) - F(\rho^K_t) - \int_{t}^{t+s} \mathcal{L}F(\rho^K_u) du \right) \prod_{i=1}^{k} h_i(\rho^K_{t_i}) } = 0,
	\end{equation}
	for all $ 0 \leq t_1 \leq \ldots \leq t_k \leq t $ and $ s \geq 0 $.
	By linearity, it is enough to consider $ F = F_{f,\phi} \in \mathcal{D}_0 $.
	Then, by Lemma~\ref{lemma:FW_finite_variance}, since $ \abs{\mathcal{L}F_{f,\phi}(\rho)} \leq C \| \phi \|^2 $,
	\begin{multline*}
		\E{ \left( F_{f,\phi}(\rho^K_{t+s}) - F_{f,\phi}(\rho^K_t) - \int_{t}^{t+s} \mathcal{L}F_{f,\phi}(\rho^K_u) du \right) \prod_{i=1}^{k} h_i(\rho^K_{t_i}) } \\ = \E{ \left( F_{f,\phi}(\rho^K_{(t+s)\wedge \tau_K}) - F_{f,\phi}(\rho^K_{t\wedge \tau_K}) - \int_{t\wedge \tau_K}^{(t+s) \wedge \tau_K} \mathcal{L}F_{f,\phi}(\rho^K_u) du \right) \prod_{i=1}^{k} h_i(\rho^K_{t_i \wedge \tau_K}) } + \littleO{1}.
	\end{multline*}
	Then, by Lemma~\ref{lemma:rescaled_generator_rho},
	\begin{multline} \label{decomp_martingale}
		\E{ \left( F_{f,\phi}(\rho^K_{(t+s)\wedge \tau_K}) - F_{f,\phi}(\rho^K_{t\wedge \tau_K}) - \int_{t\wedge \tau_K}^{(t+s) \wedge \tau_K} \mathcal{L}F_{f,\phi}(\rho^K_u) du \right) \prod_{i=1}^{k} h_i(\rho^K_{t_i \wedge \tau_K}) } \\ = \E{ \int_{t \wedge \tau_K}^{(t+s) \wedge \tau_K} \left( N_e \mathcal{S}_K(K \Nbar(u)) - 1 \right) \mathcal{L}F_{f,\phi}(\rho^K_u) du \times \prod_{i=1}^{k} h_i(\rho^K_{t_i \wedge \tau_K}) } \\ + \E{ \int_{t \wedge \tau_K}^{(t+s) \wedge \tau_K} G^K_{f,\phi}(\nu^K_{Ku}) du \times \prod_{i=1}^{k} h_i(\rho^K_{t_i \wedge \tau_K}) }.
	\end{multline}
	By \eqref{cvg_SK} and the fact that $ n \mapsto \mathcal{S}_K(K n) $ is uniformly Lipschitz on $ [n_* - \varepsilon, n_* + \varepsilon] $, we obtain that there exists a constant $ C > 0 $ such that, for all $ 0 \leq t \leq \tau_K $,
	\begin{equation} \label{bound_SK}
		\abs{ N_e \mathcal{S}_K(K \Nbar(t)) - 1 } \leq C \varepsilon, \quad \text{ almost surely.}
	\end{equation}
	In addition, for $ 0 \leq t \leq \tau_K $, using \eqref{bound_GK}, for $ K $ large enough that $ 1/K \leq \varepsilon $,
	\begin{multline*}
		\abs{ G^K_{f,\phi}(\nu^K_{K t}) } \leq \frac{4}{3} \| f^{(3)} \|_\infty \| \phi \|^3 (n_* + \varepsilon) \\ \times \left( b \sum_{k=1}^{+\infty} K^2 \left( \frac{k}{K} \right)^{2+\delta} \frac{1}{(n_*-\varepsilon)^{2+\delta}} p_k + \left( d + c(n_* + \varepsilon)  \right) \frac{1}{K (n_*-2\varepsilon)^3} \right),
	\end{multline*}
	where we have used ($ \delta \in (0,1] $ without loss of generality)
	\begin{equation*}
		\left( \frac{k}{N+k} \right)^3 \leq \left( \frac{k}{N+k} \right)^{2+\delta}.
	\end{equation*}
	By \eqref{assumption_p_finite_variance}, the sum over $ k $ on the right hand side is finite, and we obtain that there exists a constant $ C > 0 $ such that, for $ 0 \leq t \leq \tau_K $,
	\begin{equation} \label{bound_GK_almost_sure}
		\abs{G^K_{f,\phi}(\nu^K_{K t})} \leq \frac{C}{K^{\delta \wedge 1}}.
	\end{equation}
	Combining \eqref{bound_GK_almost_sure}, \eqref{bound_SK} and \eqref{decomp_martingale}, we obtain that there exists a constant $ C > 0 $ such that
	\begin{equation*}
		\limsup_{K \to \infty} \E{ \left( F(\rho^K_{t+s}) - F(\rho^K_t) - \int_{t}^{t+s} \mathcal{L}F(\rho^K_u) du \right) \prod_{i=1}^{k} h_i(\rho^K_{t_i}) } \leq C \varepsilon.
	\end{equation*}
	Since the left hand side does not depend on $ \varepsilon $, we can let it tend to zero, and we obtain \eqref{cvg_martingale_pb}, which concludes the proof of Theorem~\ref{thm:finite_variance}.
	
	\subsection{Exit time of a neighbourhood of the equilibrium population size} \label{subsec:exit_time}
	
	In this subsection, we prove Lemma~\ref{lemma:FW_finite_variance}.
	We begin by noting that, if the offspring size distribution $ (p_k, k \in \N) $ admits exponential moments, then a much stronger result follows by applying the Freidlin-Wentzell theory on the exit time of the neighbourhood of a stable equilibrium using large deviations techniques, namely that there exists $ V > 0 $ such that
	\begin{equation*}
		\lim_{K \to \infty} \P{\tau_K \leq e^{K V}} = 0.
	\end{equation*}
	This is proved for $ p_k = \delta_1 $ in \citep[Theorem~3.c]{champagnat_microscopic_2006} and more generally in \citep[Chapter~5-6]{freidlin_random_1998}.
	Since here we make a much weaker assumption on $ (p_k, k \in \N) $, we give a self-contained proof.
	
	Recall the definition of $ \tau_K $ in \eqref{def:tau_K}, and note that $ \varepsilon > 0 $ is kept fixed throughout the proof of Lemma~\ref{lemma:FW_finite_variance}.
	For convenience, we assume that $ \varepsilon < \frac{n_*}{2} $ (the probability appearing in the statement of Lemma~\ref{lemma:FW_finite_variance} is decreasing as a function of $ \varepsilon $).
	Also let $ V : (0,\infty) \to \R_+ $ be defined by
	\begin{equation} \label{def:V}
		V(n) = \frac{n}{n_*} - 1 - \log\left( \frac{n}{n_*} \right).
	\end{equation}
	(This function acts as a Lyapunov functional for the deterministic dynamical system $ \dot{n} = (bm - d - cn) n $.)
	We note that there exist constants $ \newCst{inf_V} > 0 $ and $ \newCst{sup_V} > 0 $ such that, for all $ n \in [n_* - 2\varepsilon, n_* + 2\varepsilon] $,
	\begin{equation} \label{quadratic_bound_V}
		\Cst{inf_V} (n-n_*)^2 \leq V(n) \leq \Cst{sup_V} (n-n_*)^2.
	\end{equation}
	We then set $ \newCst{sup_F} := 2 \Cst{sup_V} \varepsilon^2 $ and define for $ p \in \N $,
	\begin{equation*}
		F_p(n) := \left( V(n) \wedge \Cst{sup_F} \right)^p.
	\end{equation*}
	We can then prove the following lemma.
	
	\begin{lemma} \label{lemma:bound_E_Fp}
		Under the assumptions of Theorem~\ref{thm:finite_variance}, for all $ p \geq 1 $, there exist constants $ C_{p,1} > 0 $, $ C_{p,2} > 0 $ and $ C_{p,3} > 0 $ such that, for all $ t \geq 0 $,
		\begin{multline} \label{bound_E_Fp}
			\E{ F_p(\Nbar(t \wedge \tau_K)) } + K C_{p,1} \E{ \int_{0}^{t \wedge \tau_K} F_p(\Nbar(s)) ds } \\ \leq \E{F_p(\Nbar(0))} + C_{p,2} \E{\int_{0}^{t \wedge \tau_K} F_{p-1}(\Nbar(s)) ds} + \frac{C_{p,3}}{K^{1\wedge \delta}} t.
		\end{multline}
	\end{lemma}
	
	Before proving this lemma, let us use it to conclude the proof of Lemma~\ref{lemma:FW_finite_variance}.
	
	\begin{proof}[Proof of Lemma~\ref{lemma:FW_finite_variance}]
		Set
		\begin{equation*}
			A_p := \E{ \int_{0}^{t \wedge \tau_K} F_p(\Nbar(s)) ds }.
		\end{equation*}
		Then $ A_0 \leq t $ and, by Lemma~\ref{lemma:bound_E_Fp}, since $ F_p \geq 0 $,
		\begin{equation*}
			A_1 \leq \frac{1}{C_{1,1} K} \left( \E{F_1(\Nbar(0))} + C_{1,2} t + \frac{C_{1,3}}{K^{1\wedge \delta}} t \right).
		\end{equation*}
		Plugging this bound on the right hand side of \eqref{bound_E_Fp} with $ p = 2 $ and using the fact that $ A_2 \geq 0 $, we obtain
		\begin{multline} \label{bound_E_F2}
			\E{F_2(\Nbar(t \wedge \tau_K))} \leq \E{F_2(\Nbar(0))} \\+ \frac{C_{2,2}}{C_{1,1} K} \left( \E{F_1(\Nbar(0))} + C_{1,2} t + \frac{C_{1,3}}{K^{1\wedge \delta}} t \right) + \frac{C_{2,3}}{K^{1\wedge \delta}} t.
		\end{multline}
		We then note that, by \eqref{quadratic_bound_V},
		\begin{equation*}
			\P{\tau_K \leq t} \leq \P{F_2(\Nbar(t\wedge \tau_K)) > (\Cst{inf_V} \varepsilon^2 )^2}
		\end{equation*}
		Using Markov's inequality and \eqref{bound_E_F2}, we obtain
		\begin{multline*}
			\P{\tau_K \leq t} \leq \frac{1}{(\Cst{inf_V} \varepsilon^2)^2} \E{F_2(\Nbar(0))} \\+ \frac{C_{2,2}}{(\Cst{inf_V} \varepsilon^2)^2 C_{1,1} K} \left( \E{F_1(\Nbar(0))} + C_{1,2} t + \frac{C_{1,3}}{K^{1\wedge \delta}} t \right) + \frac{C_{2,3}}{(\Cst{inf_V} \varepsilon^2)^2 K^{1\wedge \delta}} t.
		\end{multline*}
		We conclude by noting that all the terms on the right hand side vanish as $ K \to \infty $ (using the fact that $ \Nbar(0) \to n_* $ in probability and the fact that $ F_2 $ is bounded and continuous).
	\end{proof}
	
	We complete this section by proving Lemma~\ref{lemma:bound_E_Fp}.
	
	\begin{proof}[Proof of Lemma~\ref{lemma:bound_E_Fp}]
		Note that $ (\Nbar(t), t \geq 0) $ is itself a Markov process and let $ \mathcal{B}_K $ denote its infinitesimal generator, acting on bounded real-valued functions on $ \R_+ $ as
		\begin{equation*}
			\mathcal{B}_K f(n) = K^2 b n \sum_{k=1}^{\infty} p_k \left( f\left( n + \frac{k}{K} \right) - f(n) \right) + K^2 (d + cn) n \left( f\left( n - \frac{1}{K} \right) - f(n) \right).
		\end{equation*}
		We then apply this operator to $ F_p $ and, for $ n \in [n_* - \varepsilon, n_* + \varepsilon] $, we decompose the result as follows
		\begin{multline} \label{B_K_F_p}
			\mathcal{B}_K F_p(n) = K (bm - d - cn) n F_p'(n) + \frac{1}{2} \left( b \sum_{k=1}^{\lfloor \varepsilon K \rfloor} k^2 p_k + d + cn \right) n F_p''(n) \\ + B^{(1)}_{p,K}(n) + B^{(2)}_{p,K}(n) + B^{(3)}_{p,K}(n) + B^{(4)}_{p,K}(n),
		\end{multline}
		where,
		\begin{align*}
			B^{(1)}_{p,K}(n) &:= K^2 bn \sum_{k=1}^{\lfloor \varepsilon K \rfloor} p_k \left( F_p\left( n + \frac{k}{K} \right) - F_p(n) - \frac{k}{K} F_p'(n) - \frac{k^2}{2 K^2} F_p''(n) \right) \\
			B^{(2)}_{p,K}(n) &:= K^2 (d + cn) n \left( F_p\left( n - \frac{1}{K} \right) - F_p(n) + \frac{1}{K} F_p'(n) - \frac{1}{2 K^2} F_p''(n) \right) \\
			B^{(3)}_{p,K}(n) &:= K^2 bn \sum_{k > \varepsilon K} p_k \left( F_p\left( n + \frac{k}{K} \right) - F_p(n) \right) \\
			B^{(4)}_{p,K}(n) &:= - K bn \sum_{k > \varepsilon K} k p_k F_p'(n).
		\end{align*}
		Note that, for $ n \in [n_* - \varepsilon, n_* + \varepsilon] $, by \eqref{quadratic_bound_V} and the choice of $ \Cst{sup_F} $, $ F_p(n) = V(n)^p $ and so $ F_p $ is indeed differentiable at $ n $.
		In addition $ n + \frac{k}{K} \in [n_* - 2\varepsilon, n_* + 2\varepsilon] $ for all $ -1 \leq k \leq \lfloor \varepsilon K \rfloor $, and again, by the choice of $ \Cst{sup_F} $, $ F_p\left( n + \frac{k}{K} \right) = V\left( n + \frac{k}{K} \right)^p $.
		
		Let us write $ \mathcal{V}_\varepsilon = [n_* - \varepsilon, n_* + \varepsilon] $ in what follows.
		For $ n \in \mathcal{V}_\varepsilon $,
		\begin{equation} \label{F_p_prime}
			F_p'(n) = p V'(n) V(n)^{p-1} = p \left( \frac{1}{n_*} - \frac{1}{n} \right) V(n)^{p-1}.
		\end{equation}
		Hence, using the fact that $ bm - d = c n_* $ and \eqref{quadratic_bound_V} in the second line, for $ n \in \mathcal{V}_\varepsilon $,
		\begin{align*}
			K (bm - d - cn) n F_p'(n) &= - K \frac{c p}{n_*} (n-n_*)^2 V(n)^{p-1} \\
			&\leq - K \frac{c p}{\Cst{sup_V}} V(n)^p. \numberthis \label{bound_negative_term}
		\end{align*}
		In addition, for $ n \in \mathcal{V}_\varepsilon $,
		\begin{align*}
			F_p''(n) &= p V''(n) V(n)^{p-1} + p(p-1) (V'(n))^2 V(n)^{p-2} \\
			&= \frac{p}{n^2} \left( V(n)^{p-1} + \frac{p-1}{n_*^2} (n-n_*)^2 V(n)^{p-2} \right).
		\end{align*}
		Hence, using \eqref{quadratic_bound_V}, for all $ n \in \mathcal{V}_\varepsilon $,
		\begin{equation*}
			n F_p''(n) \leq \frac{p}{(n_* - \varepsilon)^2} \left( 1 + \frac{p-1}{n_*^2 \Cst{inf_V}} \right) V(n)^{p-1}.
		\end{equation*}
		Using \eqref{assumption_p_finite_variance}, it follows that there exists a constant $ C > 0 $ (depending on $ p $ and $ \varepsilon $) such that, for all $ n \in \mathcal{V}_\varepsilon $,
		\begin{equation} \label{bound_O1_term}
			\frac{1}{2} \left( b \sum_{k=1}^{\lfloor \varepsilon K \rfloor} k^2 p_k + d + cn \right) n F_p''(n) \leq C V(n)^{p-1}.
		\end{equation}
		We then bound each $ B^{(i)}_{p,K} $ separately.
		
		\paragraph*{Bound on $ B^{(1)}_{p,K} $}
		
		By Taylor's inequality, for $ n \in \mathcal{V}_\varepsilon $ and $ 1 \leq k \leq \lfloor \varepsilon K \rfloor $,
		\begin{equation*}
			\left| F_p\left( n + \frac{k}{K} \right) - F_p(n) - \frac{k}{K} F_p'(n) - \frac{k^2}{2 K^2} F_p''(n) \right| \leq \frac{k^3}{6 K^3} \sup_{n \in \mathcal{V}_{2\varepsilon}} \left| F_p^{(3)}(n) \right|.
		\end{equation*}
		One can easily check that there exists a constant $ C > 0 $ (depending on $ p $ and $ \varepsilon $) such that
		\begin{equation*}
			\sup_{n \in \mathcal{V}_{2\varepsilon}} \left| F_p^{(3)}(n) \right| \leq C.
		\end{equation*}
		It follows that, for all $ n \in \mathcal{V}_\varepsilon $,
		\begin{equation*}
			\abs{ B^{(1)}_{p,K}(n) } \leq \frac{C b}{6} (n_* + \varepsilon) \frac{1}{K} \sum_{k=1}^{\lfloor \varepsilon K \rfloor} k^3 p_k.
		\end{equation*}
		We then use $ k^3 \leq k^{2+\delta} (\varepsilon K)^{(1-\delta)^+} $ inside the sum to obtain
		\begin{equation} \label{bound_B1}
			\abs{ B^{(1)}_{p,K}(n) } \leq \frac{C b}{6} (n_* + \varepsilon) m_{(2+\delta)} \frac{\varepsilon^{(1-\delta)^+}}{K^{1 \wedge \delta}},
		\end{equation}
		where $ m_{(2+\delta)} = \sum_{k=1}^{\infty} k^{2+\delta} p_k $.
		
		\paragraph*{Bound on $ B^{(2)}_{p,K}(n) $}
		
		Again by Taylor's inequality, for $ n \in \mathcal{V}_{\varepsilon} $,
		\begin{equation*}
			\left| F_p\left( n - \frac{1}{K} \right) - F_p(n) + \frac{1}{K} F_p'(n) - \frac{1}{2 K^2} F_p''(n) \right| \leq \frac{1}{6 K^3} \sup_{n \in \mathcal{V}_{2\varepsilon}} \left| F_p^{(3)}(n) \right|.
		\end{equation*}
		From this we deduce that, for all $ n \in \mathcal{V}_\varepsilon $,
		\begin{equation} \label{bound_B2}
			\abs{ B^{(2)}_{p,K}(n) } \leq \frac{C}{K},
		\end{equation}
		for some constant $ C > 0 $.
		
		\paragraph*{Bound on $ B^{(3)}_{p,K}(n) $}
		
		By the definition of $ F_p $, $ F_p(n) \leq \Cst{sup_F}^p $ for all $ n $, so
		\begin{equation*}
			\abs{ B^{(3)}_{p,K}(n) } \leq 2 b (n_* + \varepsilon) \Cst{sup_F}^p K^2 \P{X > \varepsilon K},
		\end{equation*}
		where $ X $ is a random variable distributed according to $ (p_k, k \in \N) $.
		By Markov's inequality,
		\begin{equation*}
			\P{X > \varepsilon K} \leq \frac{m_{(2+\delta)}}{(\varepsilon K)^{2+\delta}}.
		\end{equation*}
		As a result, for all $ n \in \mathcal{V}_\varepsilon $,
		\begin{equation} \label{bound_B3}
			\abs{ B^{(3)}_{p,K}(n) } \leq \frac{C}{K^\delta},
		\end{equation}
		for some constant $ C > 0 $ (depending on $ p $ and $ \varepsilon $).
		
		\paragraph*{Bound on $ B^{(4)}_{p,K}(n) $}
		
		First note that, by \eqref{F_p_prime}, for all $ n \in \mathcal{V}_\varepsilon $,
		\begin{equation*}
			\abs{ n F_p'(n) } \leq \frac{p}{n_*} \abs{n-n_*} V(n)^{p-1}.
		\end{equation*}
		By \eqref{quadratic_bound_V}, there exists a constant $ C > 0 $ such that $ \abs{n F_p'(n)} \leq C $ for all $ n \in \mathcal{V}_\varepsilon $.
		As a result,
		\begin{equation*}
			\abs{ B^{(4)}_{p,K}(n) } \leq C b K \E{ X \1{X > \varepsilon K} }.
		\end{equation*}
		Let $ \tilde{X} $ denote a random variable distributed according to the size-biased distribution $ (\frac{k p_k}{m}, k \in \N) $, then
		\begin{align*}
			\E{X \1{X > \varepsilon K}} &= m \P{\tilde{X} > \varepsilon K} \\
			&\leq m \frac{\E{\tilde{X}^{1+\delta}}}{(\varepsilon K)^{1+\delta}} = \frac{\E{X^{2+\delta}}}{(\varepsilon K)^{1+\delta}}.
		\end{align*}
		It follows that
		\begin{equation} \label{bound_B4}
			\abs{ B^{(4)}_{p,K}(n) } \leq \frac{C b \, m_{(2+\delta)}}{\varepsilon^{1+\delta} K^\delta}.
		\end{equation}
		
		\paragraph*{Conclusion of the proof}
		
		Combining \eqref{bound_B1}, \eqref{bound_B2}, \eqref{bound_B3} and \eqref{bound_B3}, we see that there exists $ C_{p,3} > 0 $ such that, for all $ n \in \mathcal{V}_\varepsilon $,
		\begin{equation*}
			\sum_{i=1}^{4} \abs{ B^{(i)}_{p,K}(n) } \leq \frac{C_{p,3}}{K^{1\wedge \delta}}.
		\end{equation*}
		Plugging this, along with \eqref{bound_negative_term} and \eqref{bound_O1_term}, in \eqref{B_K_F_p} we obtain that there exist constants $ C_{p,1} > 0 $ and $ C_{p,2} > 0 $ such that, for all $ n \in \mathcal{V}_\varepsilon $,
		\begin{equation*}
			\mathcal{B}_K F_p(n) \leq - K C_{p,1} V(n)^p + C_{p,2} V(n)^{p-1} + \frac{C_{p,3}}{K^{1\wedge \delta}}.
		\end{equation*}
		Moreover, by \eqref{quadratic_bound_V}, for $ n \in \mathcal{V}_\varepsilon $, $ V(n)^p = F_p(n) $, so
		\begin{equation} \label{bound_BK_Fp}
			\mathcal{B}_K F_p(n) \leq - K C_{p,1} F_p(n) + C_{p,2} F_{p-1}(n) + \frac{C_{p,3}}{K^{1\wedge \delta}}.
		\end{equation}
		Since $ F_p $ is bounded, by the optional stopping theorem,
		\begin{equation*}
			t \mapsto F_p(\Nbar(t \wedge \tau_K)) - F_p(\Nbar(0)) - \int_{0}^{t\wedge \tau_K} \mathcal{B}_K F_p(\Nbar(s)) ds
		\end{equation*}
		is a martingale and is bounded almost surely for each $ t \geq 0 $.
		Taking the expectation and plugging \eqref{bound_BK_Fp} yields \eqref{bound_E_Fp}, as desired.
	\end{proof}
	
	\section[Convergence to a Lambda-Fleming-Viot process in the alpha-stable regime]{Convergence to a $ \Lambda $-Fleming-Viot process in the $ \alpha $-stable regime} \label{sec:alpha_stable}
	
	\subsection{Stochastic averaging} \label{subsec:averaging_reminders}
	
	The proofs of Theorems~\ref{thm:alpha-stable} and \ref{thm:neveu} will be based on an application of the stochastic averaging principle as developed in \cite{kurtz_averaging_1992}.
	We recall here the main setup of \cite{kurtz_averaging_1992} and restate the main results that we will use, with some minor adaptations.
	
	Given a metric space $ S $, let $ \mathcal{M}_F(S) $ denote the space of finite measures on $ S $, and let $ \ell(S) $ denote the subspace of $ \mathcal{M}_F(\R_+ \times S) $ consisting of measures $ \mu $ such that, for all $ t \geq 0 $, $ \mu([0,t] \times S) \leq t $ (Kurtz works in the space $ \ell_m(S) $ of measures $ \mu $ for which this inequality is an equality, while we will only assume that the equality holds in the limit as $ K \to \infty $, see below).
	If $ \mu \in \ell(S) $ and $ t \geq 0 $, let $ \mu^{(t)} $  denote the restriction of $ \mu $ to $ [0,t] \times S $.
	For $ t \geq 0 $, let $ d_t $ denote the Prohorov metric on $ \mathcal{M}_F([0,t] \times S) $ and define a metric $ d $ on $ \ell(S) $ by
	\begin{equation*}
		d(\mu_1, \mu_2) := \int_{0}^{\infty} e^{-t} \left( d_t(\mu_1^{(t)}, \mu_2^{(t)}) \wedge 1 \right) dt.
	\end{equation*}
	Hence, $ \mu_n \to \mu $ in $ (\ell(S), d) $ if and only if $ \mu_n^{(t)} \to \mu^{(t)} $ for almost every $ t \geq 0 $.
	The following is proved in \cite{kurtz_averaging_1992}.
	
	\begin{lemma}[Lemma~1.3 in \cite{kurtz_averaging_1992}] \label{lemma:tightness_ell}
		A sequence of $ \ell(S) $-valued random variables $ (\Gamma_K, K \geq 1) $ is tight if and only if, for any $ \varepsilon > 0 $ and $ t \geq 0 $, there exists a compact set $ \mathcal{K} \subset S $ such that
		\begin{equation*}
			\sup_{K \geq 1} \E{ \Gamma_K([0,t] \times \mathcal{K}^c) } \leq \varepsilon.
		\end{equation*}
	\end{lemma}

	For a metric space $ E $, recall that $ C(E) $ denotes the set of continuous real-valued functions defined on $ E $, and let $ C_b(E) $ denote the set of bounded and continuous real-valued functions.
	The following is adapted from Theorem~2.1 in \cite{kurtz_averaging_1992}
	
	\begin{theorem} \label{thm:averaging}
		Let $ E_1 $ and $ E_2 $ be two complete separable metric spaces and set $ E := E_1 \times E_2 $.
		Suppose that for $ K \geq 1 $, $ \lbrace (X_K, Y_K), K \geq 1 \rbrace $ is a sequence of random variables taking values in $ D([0,\infty) \times E) $, adapted to some filtration $ (\F^K_t, t \geq 0) $.
		We assume the following.
		\begin{enumerate}
			\item \label{ass:compact_containment} The sequence $ (X_K, K \geq 1) $ satisfies the compact containment condition, i.e. for each $ \varepsilon > 0 $ and $ T > 0 $, there exists a compact set $ \mathcal{K} \subset E_1 $ such that
			\begin{equation*}
				\inf_{K \geq 1} \P{ X_K(t) \in \mathcal{K}, \forall t \in [0,T] } \geq 1-\varepsilon.
			\end{equation*}
			
			\item \label{ass:martingale_pb} There exists an operator $ A : \mathcal{D}(A) \subset C_b(E_1) \to C(E_1 \times E_2) $ and a sequence of $ \F^K_t $-stopping times $ (\tau_K, K \geq 1) $ such that, for all $ f \in \mathcal{D}(A) $, there exists a process $ (\varepsilon^f_K(t), t \geq 0) $ such that 
			\begin{equation*}
				f(X_K(t \wedge \tau_K)) - f(X_K(0)) - \int_{0}^{t\wedge \tau_K} Af(X_K(s), Y_K(s)) ds + \varepsilon^f_K(t)
			\end{equation*}
			is an $ \F^K_t $-martingale, and, for all $ T > 0 $, there exists a constant $ C_T > 0 $ such that
			\begin{equation*}
				\sup_{t \in [0, T \wedge \tau_K]} \abs{ Af(X_K(t), Y_K(t)) } \leq C_T, \quad \text{ almost surely,}
			\end{equation*}
			and a deterministic sequence $ (\eta_K, K \geq 1) $, converging to zero as $ K \to \infty $, such that
			\begin{equation*}
				\sup_{t \in [0,T]} |\varepsilon^f_K(t)| \leq \eta_K.
			\end{equation*}
		
			\item \label{ass:tightness_Y} For any $ T > 0 $, the collection of $ E_2 $-valued random variables $ \lbrace Y_K(t \wedge \tau_K), t \in [0,T], K \geq 1 \rbrace $ is tight.
		
			\item \label{ass:DA_dense} The set $ \mathcal{D}(A) $ is dense in $ C_b(E_1) $ for the topology of uniform convergence on compact sets.
			
			\item \label{ass:stopping_time} For any $ t \geq 0 $,
			\begin{equation*}
				\lim_{K \to \infty} \P{\tau_K\leq t} = 0.
			\end{equation*}
		\end{enumerate}
		We then define a sequence $ (\Gamma_K, K \geq 1) $ of $ \ell(E_2) $-valued random variables as
		\begin{equation*}
			\Gamma_K([0,t] \times B) := \int_{0}^{t \wedge \tau_K} \1{Y_K(s) \in B} ds,
		\end{equation*}
		for all measurable $ B \subset E_2 $.
		Then, under these assumptions, the sequence of $ D([0,\infty), E_1) \times \ell(E_2) $-valued random variables $ \lbrace (X_K, \Gamma_K), K \geq 1 \rbrace $ is tight and, for any limit point $ (X, \Gamma) $,
		\begin{equation} \label{Gamma_t}
			\Gamma([0,t] \times E_2) = t, \quad \forall t \geq 0,
		\end{equation}
		almost surely and there exists a filtration $ (\mathcal{G}_t, t \geq 0) $ such that, for all $ f \in \mathcal{D}(A) $,
		\begin{equation} \label{X_Gamma_mp}
			f(X(t)) - f(X(0)) - \int_{[0,t] \times E_2} Af(X(s), y) \Gamma(ds, dy)
		\end{equation}
		is a $ \mathcal{G}_t $-martingale.
	\end{theorem}

	Below, the role of $ X_K $ will be played by $ \rho_K $, and that of $ Y_K $ will be played by $ \Nbar $.
	
	The above assumptions differ slightly from those made in Theorem~2.1 of \cite{kurtz_averaging_1992}, mostly through the fact that the latter does not introduce the sequence of stopping times $ (\tau_K, K \geq 1) $ but instead assumes that
	\begin{equation*}
		\sup_{K \geq 1} \E{ \int_{0}^{T} \abs{A f(X_K(s), Y_K(s))}^p } < \infty
	\end{equation*}
	for some $ p > 1 $ and
	\begin{equation*}
		\lim_{K \to \infty} \E{ \sup_{t \in [0,T]} \abs{\varepsilon^f_K(t)} } = 0.
	\end{equation*}
	We shall see that our assumptions lead to the same conclusion as in the original result.
	The almost sure bounds on $ \varepsilon^f_K $ and $ Af(X_K(t), Y_K(t)) $ may seem restrictive, but this is accommodated by the introduction of the stopping time $ \tau_K $, which dispenses us from controlling moments of these quantities beyond the random time $ \tau_K $ (which must nonetheless diverge in probability as $ K \to \infty $).

	\begin{proof}
		Let us show briefly how to adapt the arguments of \cite{kurtz_averaging_1992} to this slightly different setting.
		Let us set $ \tilde{X}_K(t) := X_K(t \wedge \tau_K) $.
		By assumptions~\ref{ass:compact_containment} and \ref{ass:martingale_pb} and by Theorem~3.9.1 in \cite{ethier_markov_1986}, $ (\tilde{X}_K, K \geq 1) $ is tight in $ D([0,\infty) \times E_1) $.
		Then, by assumption~\ref{ass:stopping_time}, $ (X_K, K \geq 1) $ is tight and any limit point is also a limit point of $ (\tilde{X}_K, K \geqq 1) $, and vice versa.
		
		By assumption~\ref{ass:tightness_Y}, for any $ \varepsilon > 0 $ and $ T \geq 0 $, there exists a compact set $ \mathcal{K} \subset E_2 $ such that
		\begin{equation*}
			\sup_{K \geq 1} \sup_{t \in [0,T]} \P{ Y_K(t) \notin \mathcal{K} } \leq \varepsilon.
		\end{equation*}
		It follows that
		\begin{equation*}
			\sup_{K \geq 1} \E{ \Gamma_K([0,T] \times \mathcal{K}^c) } \leq \varepsilon t.
		\end{equation*}
		Hence, by Lemma~\ref{lemma:tightness_ell}, $ (\Gamma_K, K \geq 1) $ is tight in $ \ell(E_2) $.
		
		Let $ (X, \Gamma) $ be a limit point of $ \lbrace (X_K, \Gamma_K), K \geq 1 \rbrace $.
		Assumption~\ref{ass:stopping_time} directly entails \eqref{Gamma_t}.
		Define the filtration $ (\mathcal{G}_t, t \geq 0) $ as
		\begin{equation*}
			\mathcal{G}_t := \sigma\lbrace X(s), \Gamma([0,s] \times B), s \in [0,t], B \in \mathcal{B}(E_2) \rbrace.
		\end{equation*}
		For $ f \in \mathcal{D}(A) $, set
		\begin{equation*}
			Z_K(t) := f(X_K(t\wedge \tau_K)) - f(X_K(0)) - \int_{[0,t] \times E_2} Af(X_K(s), y) \Gamma_K(ds, dy),
		\end{equation*}
		and
		\begin{equation*}
			Z(t) := f(X(t)) - f(X(0)) - \int_{[0,t] \times E_2} Af(X(s), y) \Gamma(ds, dy).
		\end{equation*}
		Then, by Lemma~1.5 in \cite{kurtz_averaging_1992}, $ Z_K \to Z $ in distribution in $ D([0,\infty), \R) $ along an appropriate subsequence.
		Moreover, by assumption~\ref{ass:martingale_pb}, $ Z_K + \varepsilon^f_K $ is an $ \F^K_t $-martingale and $ Z_K $ is uniformly integrable (even bounded almost surely).
		Since $ \varepsilon^f_K \to 0 $ in $ L^1 $, it follows that $ Z $ is a $ \mathcal{G}_t $-martingale.
		This concludes the proof of Theorem~\ref{thm:averaging}.
	\end{proof}
	
	We complement this result by the following, which is adapted from Example~2.3 in \cite{kurtz_averaging_1992}.
	
	\begin{proposition} \label{prop:averaging_Y}
		Assume that, in addition to the assumptions of Theorem~\ref{thm:averaging}, there exist an operator $ \mathcal{B} : \mathcal{D}(\mathcal{B}) \subset C_b(E_2) \to C_b(E_1, E_2) $ and a sequence $ (\beta_K, K \geq 1) $ such that $ \beta_K \to \infty $ as $ K \to \infty $ and, for all $ g \in \mathcal{D}(\mathcal{B}) $, there exists a process $ (\delta^g_K(t), t \geq 0) $ such that
		\begin{equation*}
			g(Y_K(t \wedge \tau_K)) - g(Y_K(0)) - \beta_K \int_{0}^{t \wedge \tau_K} \mathcal{B}g(X_K(s), Y_K(s)) ds + \delta^g_K(t)
		\end{equation*}
		is an $ \F^K_t $-martingale and, for all $ T > 0 $, there exists a deterministic sequence $ (\eta_K, K \geq 1) $ such that $ \eta_K = \littleO{\beta_K} $ and
		\begin{equation*}
			\sup_{t \in [0, T]} | \delta^g_K(t) | \leq \eta_K,
		\end{equation*}
		almost surely.
		Then, for all $ g \in \mathcal{D}(\mathcal{B}) $, for any limit point $ (X, \Gamma) $ of the sequence $ \lbrace (X_K, \Gamma_K), K \geq 1 \rbrace $,
		\begin{equation*}
			\int_{[0,t] \times E_2} \mathcal{B}g(X(s), y) \Gamma(ds, dy) = 0, \quad \forall t \geq 0,
		\end{equation*}
		almost surely.
	\end{proposition}
	
	The above proposition follows by exactly the same arguments as in Example~2.3 of \cite{kurtz_averaging_1992}.
	It will be used to characterise the limit points of $ (\Gamma_K, K \geq 1) $.
	We note also that, in our applications below, the dynamics of $ Y_K $ (a.k.a. $ \Nbar $) do not depend on $ X_K $ (a.k.a. $ \rho^K $), so that $ \mathcal{B}g(x, y) = \mathcal{B}g(y) $.
	Applying Proposition~\ref{prop:averaging_Y}, we shall thus conclude that the limit $ \Gamma $ is deterministic and concentrated on the equilibrium point $ n_* $, i.e. $ \Gamma(ds, dy) = \delta_{n_*}(dy) ds $ almost surely.
	Plugging this in \eqref{X_Gamma_mp} will allow us to identify the limit points of the sequence $ (X_K, K \geq 1) $ as solutions to a classical martingale problem.
	
	\subsection{Proof of the main result} \label{subsec:proof_stable}
	
	Recall that, in the case of Theorem~\ref{thm:alpha-stable}, for $ t \geq 0 $,
	\begin{align*}
		\Nbar(t) := \frac{1}{K} \langle \nu^K_{K^{\alpha-1} t}, 1 \rangle, && \rho^K_t := \varrho(\nu^K_{K^{\alpha-1} t}).
	\end{align*}
	Let $ (\F^K_t, t \geq 0) $ denote the natural filtration associated to the process $ \lbrace (\rho^K_t, \Nbar(t)), t \geq 0 \rbrace $.
	
	Note that, since $ (\Nbar(0), K \geq 1) $ is tight in $ (0,\infty) $, for any $ \varepsilon $, there exist two constants $ c_0 > 0 $ and $ C_0 > 0 $ and another sequence of processes $ \lbrace (\rho^{K,\varepsilon}_t, \overline{N}_{K,\varepsilon}(t)), t \geq 0 \rbrace $ satisfying all the assumptions of Theorem~\ref{thm:alpha-stable} and such that
	\begin{enumerate}
		\item $ c_0 \leq \overline{N}_{K,\varepsilon}(0) \leq C_0 $ for all $ K \geq 1 $, almost surely,
		\item for all $ K \geq 1 $, $ \P{ (\rho^K, \Nbar) = (\rho^{K,\varepsilon}, \overline{N}_{K,\varepsilon}) } \geq 1-\varepsilon $.
	\end{enumerate}
	If we then show that the statement of Theorem~\ref{thm:alpha-stable} holds for $ \lbrace (\rho^{K,\varepsilon}, \overline{N}_{K,\varepsilon}), K \geq 1 \rbrace $ for any $ \varepsilon $, the statement will also hold for $ \lbrace (\rho^K, \Nbar), K \geq 1 \rbrace $.
	As a result, without loss of generality, we assume that there exist $ c_0 \in (0, n_*) $ and $ C_0 > 0 $ such that, for all $ K \geq 1 $,
	\begin{equation} \label{N0_compact}
		c_0 \leq \Nbar(0) \leq C_0,
	\end{equation}
	almost surely.
	
	We start by stating two preliminary lemmas on the rescaled population size process $ (\Nbar(t), t \geq 0) $.
	
	\begin{lemma} \label{lemma:bounded_expectation}
		Under the assumptions of Theorem~\ref{thm:alpha-stable} and \eqref{N0_compact}, there exists a constant $ \newCst{Exp_Nbar} > 0 $ such that, for all $ K \geq 1 $,
		\begin{equation*}
			\sup_{t \geq 0} \E{\Nbar(t)} \leq \Cst{Exp_Nbar}.
		\end{equation*}
	\end{lemma}
	
	We also define a stopping time $ \tau_K $ as
	\begin{equation} \label{def:tau_K_alpha}
		\tau_K := \inf \left\lbrace t \geq 0 : \Nbar(t) < \frac{c_0}{2} \right\rbrace.
	\end{equation}
	We then have the following.
	
	\begin{lemma} \label{lemma:tau_K_stable}
		Under the assumptions of Theorem~\ref{thm:alpha-stable} and \eqref{N0_compact}, for any $ t \geq 0 $,
		\begin{equation*}
			\lim_{K \to \infty} \P{\tau_K \leq t} = 0.
		\end{equation*}
	\end{lemma}
	
	Lemmas~\ref{lemma:bounded_expectation} and \ref{lemma:tau_K_stable} are proved in Subsection~\ref{subsec:first_moment} below.
	
	~~
	
	Define a linear operator $ A : \mathcal{D}_0 \subset C(\M) \to C(\M, \R_+) $ by
	\begin{equation} \label{def:A}
		A F_{f,\phi}(\rho, n) := n b p_0 \int_{0}^{\infty} \int_\X \left( f\left( \frac{n \langle \rho, \phi \rangle + y \phi(x)}{n + y} \right) - f(\langle \rho, \phi \rangle) \right) \rho(dx) \frac{dy}{y^{1+\alpha}},
	\end{equation}
	for all $ F_{f,\phi} \in \mathcal{D}_0 $, and extend $ A $ to $ \mathcal{D} $ by linearity.
	By a change of variables, the above is also
	\begin{equation} \label{A_rewritten}
		A F_{f, \phi}(\rho, n) = n^{1-\alpha} b p_0 \int_{0}^{1} \int_\X \left( f\left( (1-u) \langle \rho, \phi \rangle + u \phi(x) \right) - f(\langle \rho, \phi \rangle) \right) \rho(dx) (1-u)^{1-\alpha} u^{-1-\alpha} du,
	\end{equation}
	where we recognise, when $ n $ is fixed, the generator of the Fleming-Viot process which is dual to the $ Beta(\alpha, 2-\alpha) $ coalescent.
	We then have the following lemma, which we prove in Subsection~\ref{subsec:averaging}.
	
	\begin{lemma} \label{lemma:averaging_rho}
		Under the assumptions of Theorem~\ref{thm:alpha-stable} and \eqref{N0_compact}, for any $ F = F_{f,\phi} \in \mathcal{D}_0 $, there is a process $ (\varepsilon^{f,\phi}_K(t), t \geq 0) $ such that
		\begin{equation*}
			F_{f,\phi}(\rho^K_{t \wedge \tau_K}) - F_{f, \phi}(\rho^K_0) - \int_{0}^{t \wedge \tau_K} A F_{f,\phi}(\rho^K_s, \Nbar(s)) ds - \varepsilon^{f, \phi}_K(t)
		\end{equation*}
		is an $ \F^K_t $-martingale.
		In addition, for any $ T > 0 $, there exists a deterministic sequence $ (\eta_K, K \geq 1) $, converging to zero as $ K \to \infty $ and such that, for all $ K \geq 1 $,
		\begin{equation*}
			\sup_{t \in [0, T]} \abs{\varepsilon^{f,\phi}_K(t)} \leq \eta_K, \quad \text{ almost surely.}
		\end{equation*}
	\end{lemma}
	
	Theorem~\ref{thm:alpha-stable} will follow from the above result, combined with the following lemma, which allows us to show that the occupation measure of $ (\Nbar(t), t \geq 0) $ concentrates on $ \lbrace n_* \rbrace $ as $ K \to \infty $.
	Let $ \mathcal{B} : C^1(\R_+) \to C(\R_+) $ be the operator defined as
	\begin{equation} \label{def:B}
		\mathcal{B} g(n) := (b m - d - cn) n g'(n).
	\end{equation}
	In words, $ \mathcal{B} $ is the generator of the deterministic process which follows the flow of the logistic equation on $ \R_+ $.
	For $ \lambda > 0 $, let $ g_{\lambda} : \R_+ \to [0,1] $ be the function defined by
	\begin{align*}
		g_{\lambda}(n) := e^{-\lambda n}.
	\end{align*}
	
	\begin{lemma} \label{lemma:Nbar}
		Under the assumptions of Theorem~\ref{thm:alpha-stable}, for any $ \lambda > 0 $, there exists a process $ (\delta^\lambda_K(t), t \geq 0) $ such that
		\begin{equation*}
			g_\lambda(\Nbar(t)) - g_\lambda(\Nbar(0)) - K^{\alpha-1} \int_{0}^{t} \mathcal{B} g_\lambda(\Nbar(s)) ds + \delta^\lambda_K(t)
		\end{equation*}
		is an $ \F^K_t $-martingale.
		In addition, for any $ T > 0 $, there exists a constant $ C_T > 0 $ such that, for all $ K \geq 1 $,
		\begin{equation*}
			\sup_{t \in [0, T]} \abs{\delta^\lambda_K(t)} \leq C_T, \quad \text{ almost surely.}
		\end{equation*}
	\end{lemma}
	
	Lemmas~\ref{lemma:averaging_rho} and \ref{lemma:Nbar} are proved in Subsections~\ref{subsec:averaging} and \ref{subsec:Nbar}, respectively.
	Let us now show how Theorem~\ref{thm:alpha-stable} follows from the above lemmas.
	
	\begin{proof}[Proof of Theorem~\ref{thm:alpha-stable}]
		Theorem~\ref{thm:alpha-stable} then follows by the above lemmas and Theorem~\ref{thm:averaging} and Proposition~\ref{prop:averaging_Y} above, taking $ X_K(t) = \rho^K_t $ (with $ E_1 = \M $) and $ Y_K(t) = \Nbar(s) $.
		Let $ \Gamma_K $ denote the random measure on $ \R_+ \times \R_+ $ defined by
		\begin{equation*}
			\Gamma_K([0,t] \times B) := \int_{0}^{t \wedge \tau_K} \1{\Nbar(s) \in B} ds,
		\end{equation*}
		for all measurable $ B \subset \R_+ $ and $ t \geq 0 $.
		Then $ \Gamma_K $ is a random variable taking values in $ \ell(\R_+) $.
		
		Assumption~\ref{ass:compact_containment} of Theorem~\ref{thm:averaging} is then automatically satisfied since $ \M $ is compact.
		Assumption~\ref{ass:martingale_pb} follows from Lemma~\ref{lemma:averaging_rho} and by linearity, also noting that by \eqref{A_rewritten}, for any $ F_{f,\phi} \in \mathcal{D}_0 $, there exists a constant $ C > 0 $ such that, for all $ n \geq 0 $ and $ \rho \in \M $,
		\begin{equation*}
			\abs{A F_{f,\phi}(\rho, n)} \leq \frac{C}{n^{\alpha-1}}.
		\end{equation*}
		Assumption~\ref{ass:tightness_Y} is satisfied by Lemma~\ref{lemma:bounded_expectation}.
		Finally, $ \mathcal{D} $ is dense in $ C(\M) $ and Assumption~\ref{ass:stopping_time} follows from Lemma~\ref{lemma:tau_K_stable}.
		
		By Theorem~\ref{thm:averaging}, $ \lbrace (\rho^K, \Gamma_K), K \geq 1 \rbrace $ is tight in $ D([0,\infty), \M) \times \ell(\R_+) $, and, for any limit point $ (\rho, \Gamma) $ and for all $ F \in \mathcal{D} $,
		\begin{equation} \label{martingale_rho}
			F(\rho_t) - F(\rho_0) - \int_{[0,t] \times \R_+} A F(\rho_s, n) \Gamma(ds, dn)
		\end{equation}
		is an $ \F_t $-martingale, where
		\begin{equation*}
			\F_t := \sigma \lbrace \rho_s, \Gamma([0,s] \times B) : s \leq t, B \in \mathcal{B}(\R_+) \rbrace.
		\end{equation*}
		In addition, by Lemma~\ref{lemma:Nbar_neveu}, we can apply Proposition~\ref{prop:averaging_Y} with $ \beta_K = K^{\alpha-1} $ to obtain that, for any $ \lambda > 0 $,
		\begin{equation*}
			\int_{[0,t] \times \R_+} \left( b m - d - c n \right) n e^{-\lambda n} \Gamma(ds, dn) = 0, \quad \text{ for all } t \geq 0,
		\end{equation*}
		almost surely.
		It follows that $ \Gamma([0,t] \times dn) $ is concentrated on $ \lbrace 0, n_* \rbrace $.
		In addition, by the definition of $ \Gamma_K $ and $ \tau_K $, for any continuous function $ f : \R_+ \to \R $ supported on $ [0, \frac{c_0}{2}) $,
		\begin{equation*}
				\int_{[0,t] \times \R_+} f(n) \Gamma(ds, dn) = 0, \quad \text{ for all } t \geq 0,
			\end{equation*}
		almost surely.
		Combined with \eqref{Gamma_t}, this shows that
		\begin{equation*}
			\Gamma(ds, dn) = \delta_{n_*}(dn) ds,
		\end{equation*}
		almost surely, which proves \eqref{bound_pop_size_alpha-stable}.
		Plugging this in \eqref{martingale_rho}, we obtain that $ (\rho_t, t \geq 0) $ solves the martingale problem associated to $ (A, \mathcal{D}) $, which concludes the proof of Theorem~\ref{thm:alpha-stable}.
	\end{proof}
	
	\subsection{First moment bound and control of the stopping time} \label{subsec:first_moment}
	
	Let $ \mathcal{B}_K $ denote the infinitesimal generator of $ (\Nbar(t), t \geq 0) $, i.e.
	\begin{equation} \label{def:BK}
		\mathcal{B}_K g(n) = K^\alpha bn \sum_{k=1}^{\infty} p_k \left( g\left( n + \frac{k}{K} \right) - g(n) \right) + K^\alpha (d + cn) n \left( g\left( n - \frac{1}{K} \right) - g(n) \right).
	\end{equation}
	For $ k_0 \in \N \setminus \lbrace 0 \rbrace $, let $ \mathcal{B}_K^{(k_0)} $ denote the corresponding generator with no birth of more than $ k_0 $ individuals, i.e.
	\begin{equation*}
		\mathcal{B}_K^{(k_0)} g(n) = K^\alpha bn \sum_{k=1}^{k_0} p_k \left( g\left( n + \frac{k}{K} \right) - g(n) \right) + K^\alpha (d + cn) n \left( g\left( n - \frac{1}{K} \right) - g(n) \right).
	\end{equation*}
	We then have the following.
	
	\begin{proposition} \label{prop:coupling}
		For any fixed $ K \geq 1 $, let $ (N_{k}, k \geq 1) $ be a family of random variables such that $ N_{k} $ takes values in $ \frac{1}{K} \N $ and, for all $ k \geq 1 $,
		\begin{equation*}
			N_{k} \leq N_{k+1} \leq \ldots \leq \Nbar(0),
		\end{equation*}
		almost surely.
		Then there exists a sequence of processes $ \lbrace(\Nbar^{(k)}(t), t \geq 0), k \geq 1 \rbrace $ such that, for each $ k \geq 1 $, $ (\Nbar^{(k)}(t), t \geq 0) $ is a Markov process with generator $ \mathcal{B}^{(k)}_K $, started from $ N_{k} $ and, for all $ k \geq 1 $
		\begin{equation*}
			\Nbar^{(k)}(t) \leq \Nbar^{(k+1)}(t) \leq \ldots \leq \Nbar(t), \quad \forall t \geq 0,
		\end{equation*}
		almost surely.
		If moreover $ N_{k} \to \Nbar(0) $ as $ k \to \infty $ almost surely, then the sequence can be constructed such that $ \Nbar^{(k)}(\cdot) \to \Nbar(\cdot) $ locally uniformly as $ k \to \infty $, almost surely.
	\end{proposition}
	
	The proof of this proposition is a classical coupling result and is thus omitted.
	We now prove Lemma~\ref{lemma:bounded_expectation}.
	
	\begin{proof}[Proof of Lemma~\ref{lemma:bounded_expectation}]
		We would like to use the fact that 
		\begin{equation} \label{martingale_Nbar}
			M_K(t) := \Nbar(t) - \Nbar(0) - \int_{0}^{t} \left( bm - d - c \Nbar(s) \right) \Nbar(s) ds
		\end{equation}
		is a local martingale.
		However, we need to be careful when taking expectations of this quantity, as $ \Nbar(t) $ does not have a finite second moment.
		To circumvent this, fix $ K \geq 1 $ and let $ (\Nbar^{(k)}(t), t \geq 0) $ be given as in Proposition~\ref{prop:coupling} for a sequence $ N_{k} $ such that $ N_{k} \to \Nbar(0) $ as $ k \to \infty $, almost surely, so that $ \Nbar^{(k)}(t) \to \Nbar(t) $ as $ k \to \infty $ almost surely for all $ t \geq 0 $.

		Then, for any $ k \geq 1 $, $ \Nbar^{(k)}(t) $ has finite moments of all orders.
		Let $ m^{(k)}(t) $ denote its expectation, i.e.
		\begin{equation*}
			m^{(k)}(t) := \E{ \Nbar^{(k)}(t) },
		\end{equation*}
		and set
		\begin{equation*}
			m_{(k)} := \sum_{\ell=1}^{k} \ell p_\ell.
		\end{equation*}
		By the definition of $ (\Nbar^{(k)}(t), t \geq 0) $,
		\begin{equation*}
			Z^{(k)}(t) := \Nbar^{(k)}(t) - N_k - K^{\alpha-1} \int_{0}^{t} \left( b m_{(k)} - d - c \Nbar^{(k)}(s) \right) \Nbar^{(k)}(s) ds
		\end{equation*}
		is a square integrable martingale.
		Taking the expectation of $ Z^{(k)}(t) - Z^{(k)}(s) $ and using Jensen's inequality yields
		\begin{equation*}
			m^{(k)}(t) \leq m^{(k)}(s) + K^{\alpha-1} \int_{s}^{t} (b m_{(k)} - d - c m^{(k)}(u)) m^{k}(u) du.
		\end{equation*}
		It follows that $ m^{(k)}(t) \leq m^{(k)}_+(K^{\alpha-1} t) $, where
		\begin{equation*}
			\left\lbrace
			\begin{aligned}
				&\deriv{m^{(k)}_+(t)}{t} = (b m_{(k)} - d - c m^{(k)}_+(t)) m^{(k)}_+(t), \\
				& m^{(k)}_+(0) = \E{N_k}.
			\end{aligned}
			\right.
		\end{equation*}
		We then observe that
		\begin{equation*}
			\sup_{t \geq 0} m^{(k)}_+(t) \leq \E{N_k} \vee \left( \frac{b m_{(k)}-d}{c} \right).
		\end{equation*}
		On the one hand, since $ m = \sum_{k=1}^{\infty} k p_k < \infty $, $ m_{(k)} \to m $ as $ k \to \infty $, and, on the other hand, by the monotone convergence theorem, $ \E{N_k} \to \E{\Nbar(0)} $ as $ k \to \infty $ (moreover, both sequences are non-decreasing), and $ \sup_{K \geq 1} \E{\Nbar(0)} < \infty $ by \eqref{N0_compact}.
		It follows that, for all $ k \geq 1 $,
		\begin{equation} \label{bound_m_k}
			\sup_{t \geq 0} m^{(k)}(t) \leq C_0 \vee n_*.
		\end{equation}
		Using the monotone convergence theorem again, we obtain that $ m^{(k)}(t) \to \E{\Nbar(t)} $ as $ k \to \infty $ for all $ t \geq 0 $.
		By \eqref{bound_m_k}, we obtain that, for all $ t \geq 0 $,
		\begin{equation*}
			\sup_{K \geq 1} \E{\Nbar(t)} \leq C_0 \vee n_*,
		\end{equation*}
		yielding the result.
	\end{proof}

	We now prove Lemma~\ref{lemma:tau_K_stable}.
	
	\begin{proof}[Proof of Lemma~\ref{lemma:tau_K_stable}]
		Recall that $ c_0 < n_* $.
		We then take $ k_0 \in \N $ such that
		\begin{equation*}
			n_0 := \frac{1}{c} \left( b \sum_{k=1}^{k_0} k p_k - d \right) \in \left( \frac{c_0}{2}, c_0 \right).
		\end{equation*}
		Then let $ (\Nbar^{(k_0)}(t), t \geq 0) $ be given by Proposition~\ref{prop:coupling}, with an initial condition $ \Nbar^{(k_0)}(0) $ such that
		\begin{equation} \label{cvg_N_bar_k0_0}
			\Nbar^{(k_0)}(0) \cvgas{K} n_0,
		\end{equation}
		in probability, and such that $ \Nbar^{(k_0)}(t) \leq \Nbar(t) $ for all $ t \geq 0 $ almost surely (this is possible by \eqref{N0_compact} and the fact that $ n_0 < c_0 $).
		By Lemma~\ref{lemma:FW_finite_variance}, we know that, for any $ \varepsilon > 0 $ and $ T > 0 $,
		\begin{equation*}
			\lim_{K \to \infty} \P{ \sup_{t \in [0,KT]} \abs{ \Nbar^{(k_0)}(t) - n_0} > \varepsilon } = 0.
		\end{equation*}
		Choosing $ \varepsilon $ small enough that $ n_0 - \varepsilon > \frac{c_0}{2} $, we obtain that
		\begin{equation} \label{lower_bound}
			\lim_{K \to \infty} \P{ \Nbar(t \wedge \tau_K) < \frac{c_0}{2} } = 0,
		\end{equation}
	which concludes the proof of the lemma.
	\end{proof}
	
	\subsection{Stochastic averaging for the distribution of types} \label{subsec:averaging}
	
	In this subsection, we prove Lemma~\ref{lemma:averaging_rho}.
	For this purpose, we introduce the following notation.
	For $ F_{f,\phi} \in \mathcal{D}_0 $, $ \rho \in \M $, $ n > 0 $ and $ y \geq -1 $ such that $ n + y > 0 $,
	\begin{equation} \label{def:Theta}
		\Theta_{n,y} F_{f,\phi}(\rho) := \int_\X \left( f\left( \frac{n\langle \rho, \phi \rangle + y \phi(x)}{n + y} \right) - f(\langle \rho, \phi \rangle) \right) \rho(dx).
	\end{equation}
	Then define $ A_K : \mathcal{D}_0 \to C(\M \times \R_+) $ as
	\begin{equation*}
		A_K F_{f,\phi}(\rho, n) := b n K^\alpha \sum_{k=1}^{\infty} p_k \Theta_{n, k/K} F_{f,\phi}(\rho) + (d + cn) n K^\alpha \Theta_{n, -1/K} F_{f,\phi}(\rho).
	\end{equation*}
	We then see that, by definition,
	\begin{equation*}
		F_{f,\phi}(\rho^K_t) - F_{f,\phi}(\rho^K_0) - \int_{0}^{t} A_K F_{f,\phi}(\rho^K_s, \Nbar(s)) ds
	\end{equation*}
	is a local martingale with respect to the filtration $ (\F^K_t, t \geq 0) $.
	
	Let us also note that, if we define, for $ y \geq 1 $,
	\begin{equation} \label{def:p_y}
		p(y) := y^{1+\alpha} p_{\lfloor y \rfloor},
	\end{equation}
	then, for any bounded function $ f : \N \to \R $,
	\begin{equation*}
		\sum_{k=1}^{\infty} f(k) p_k = \int_{1}^{\infty} f(\lfloor y \rfloor) \frac{p(y)}{y^{1+\alpha}} dy.
	\end{equation*}
	Moreover, by the assumptions of Theorem~\ref{thm:alpha-stable},
	\begin{equation} \label{cvg_p_y}
		\lim_{y \to \infty} p(y) = p_0,
	\end{equation}
	and $ p $ is bounded on $ [1,\infty) $.
	This notation will be used throughout this section and the next to simplify many computations.
	Let also $ C_p > 0 $ be defined as
	\begin{equation} \label{def:Cp}
		C_p := \sup_{y \in [1,\infty)} p(y).
	\end{equation}
	In particular, by a change of variables,
	\begin{equation} \label{integral_representation_sum_pk}
		K^\alpha \sum_{k=1}^{\infty} f\left( \frac{k}{K} \right) p_k = \int_{\frac{1}{K}}^{\infty} f\left( y_K \right) \frac{p(K y)}{y^{1+\alpha}} dy,
	\end{equation}
	where $ y_K = \frac{\lfloor K y \rfloor}{K} $.
	
	We can now prove Lemma~\ref{lemma:averaging_rho}.
	
	\begin{proof}[Proof of Lemma~\ref{lemma:averaging_rho}]
		Using the notation introduced above and \eqref{integral_representation_sum_pk}, we can write $ A_K F_{f,\phi} $ as
		\begin{equation*}
			A_K F_{f,\phi}(\rho, n) = b n \int_{\frac{1}{K}}^{\infty} \Theta_{n, y_K} F_{f,\phi}(\rho) \frac{p(Ky)}{y^{1+\alpha}} dy + (d + cn) n K^\alpha \Theta_{n, -1/K} F_{f,\phi}(\rho).
		\end{equation*}
		We also note that, by the definition of $ A $ in \eqref{def:A},
		\begin{equation*}
			A F_{f,\phi}(\rho, n) = b n \int_{0}^{\infty} \Theta_{n,y} F_{f,\phi}(\rho) \frac{p_0}{y^{1+\alpha}} dy.
		\end{equation*}
		Subtracting this to the previous equation, we obtain
		\begin{equation} \label{diff_AK_A}
			A_K F_{f,\phi}(\rho, n) - A F_{f,\phi}(\rho, n) = \sum_{i=1}^{4} R^{f,\phi}_{K,i}(\rho, n),
		\end{equation}
		where
		\begin{align*}
			R^{f,\phi}_{K,1}(\rho, n) &:= - b n p_0 \int_{0}^{\frac{1}{K}} \Theta_{n,y} F_{f,\phi}(\rho) \frac{dy}{y^{1+\alpha}}, \\
			R^{f,\phi}_{K,2}(\rho, n) &:= b n \int_{\frac{1}{K}}^{\infty} \Theta_{n, y_K} F_{f,\phi}(\rho) \frac{p(Ky) - p_0}{y^{1+\alpha}} dy, \\
			R^{f,\phi}_{K,3}(\rho, n) &:= b n p_0 \int_{\frac{1}{K}}^{\infty} \left( \Theta_{n, y_K} F_{f, \phi}(\rho) - \Theta_{n,y} F_{f,\phi}(\rho) \right) \frac{dy}{y^{1+\alpha}}, \\
			R^{f,\phi}_{K,4}(\rho, n) &:= (d + cn) n K^\alpha \Theta_{n, -1/K} F_{f,\phi}(\rho).
		\end{align*}
		We now bound each term separately, uniformly for $ n \geq c_0 / 2 $.
		
		\paragraph*{Bound on $ R^{f,\phi}_{K,1}(\rho, n) $}
		
		We start by noting that, for $ F_{f,\phi} \in \mathcal{D} $,
		\begin{equation} \label{bound_Theta_small_y}
			\abs{\Theta_{n,y} F_{f,\phi}(\rho) } \leq \| f'' \| \| \phi \|^2 \left( \frac{y}{n + y} \right)^2.
		\end{equation}
		It follows that
		\begin{align*}
			\abs{ R^{f,\phi}_{K,1}(\rho, n) } &\leq \frac{b p_0}{n} \| f'' \| \| \phi \|^2 \int_{0}^{\frac{1}{K}} y^{1-\alpha} dy \\
			&= \frac{b p_0}{n (2-\alpha)} \| f'' \| \| \phi \|^2 \frac{1}{K^{2-\alpha}}.
		\end{align*}
		As a result, for all $ n \geq c_0 / 2 $,
		\begin{equation} \label{bound_R_K1}
			\abs{ R^{f,\phi}_{K,1}(\rho, n) } \leq \frac{2 b p_0}{c_0 (2-\alpha)} \| f'' \| \| \phi \|^2 \frac{1}{K^{2-\alpha}}.
		\end{equation}
		
		\paragraph*{Bound on $ R^{f,\phi}_{K,2}(\rho, n) $}
		
		Using \eqref{bound_Theta_small_y} again and splitting the integral over $ y $, we obtain
		\begin{equation*}
			\abs{R^{f,\phi}_{K,2}(\rho, n)} \leq \frac{b}{n} \| f'' \| \| \phi \|^2 \int_{1/K}^{1} y^{1-\alpha} \abs{p(K y) - p_0} dy + b \| f'' \| \| \phi \|^2 \int_{1}^{\infty} \abs{p(Ky) - p_0} \frac{dy}{y^\alpha}.
		\end{equation*}
		We then see that both terms on the right hand side tend to zero as $ K \to \infty $ by \eqref{cvg_p_y} and the dominated convergence theorem.
		Moreover, this convergence is uniform over $ n \geq c_0 / 2 $, so there exists a deterministic sequence $ (\eta_K, K \geq 1) $ converging to zero as $ K \to \infty $ such that, for all $ n \geq c_0 / 2 $,
		\begin{equation} \label{bound_R2}
			\abs{ R^{f,\phi}_{K,2}(\rho, n) } \leq \eta_K.
		\end{equation}
		
		\paragraph*{Bound on $ R^{f,\phi}_{K,3}(\rho, n) $}
		
		For this bound, we note that
		\begin{multline} \label{bound_f_yK_y}
			\abs{ f\left( \frac{n \langle \rho, \phi \rangle + y_K \phi(x)}{n + y_K} \right) - f\left( \frac{n \langle \rho, \phi \rangle + y \phi(x)}{n + y} \right) } \\
			\begin{aligned}
				&\leq \| f' \| \left( \abs{ \frac{n}{n + y_K} - \frac{n}{n + y} } \| \phi \| + \abs{ \frac{y_K}{n + y_K} - \frac{y}{n + y} } \| \phi \| \right) \\
				&= 2 \| f' \| \| \phi \| n \frac{\abs{y_K - y}}{(n+ y) (n+y_K)} \\
				&\leq \frac{2}{n} \| f' \| \| \phi \| \abs{y_K - y}.
			\end{aligned}
		\end{multline}
		Furthermore, for all $ y \geq 0 $, $ \abs{y_K - y} \leq \frac{1}{K} $.
		This yields
		\begin{equation} \label{bound_R3_large_y}
			b p_0 n \int_{1}^{\infty} \abs{ \Theta_{n, y_K} F_{f,\phi}(\rho) - \Theta_{n,y} F_{f,\phi}(\rho) } \frac{dy}{y^{1+\alpha}} \leq 2 b p_0 \| f' \| \| \phi \| \frac{1}{K} \int_{1}^{\infty} \frac{dy}{y^{1+\alpha}}.
		\end{equation}
		To bound the integral for $ y \leq 1 $, we need a sharper bound.
		First note that, by Taylor's inequality, for $ u, u' \in [0,1] $,
		\begin{multline*}
			| f((1-u') \langle \rho, \phi \rangle + u'\phi(x)) - f((1-u) \langle \rho, \phi \rangle + u \phi(x)) \\ - (u'-u) (\phi(x) - \langle \rho, \phi \rangle) f'((1-u)\langle \rho, \phi \rangle + u \phi(x)) | \\ \leq \frac{1}{2} \| f'' \| (u'-u)^2 (\phi(x) - \langle  \rho, \phi \rangle)^2.
		\end{multline*}
		In addition, using Taylor's inequality again,
		\begin{equation*}
			| f'((1-u) \langle \rho, \phi \rangle + u \phi(x)) - f'(\langle \rho, \phi \rangle) | \leq \| f'' \| u | \phi(x) - \langle \rho, \phi \rangle |.
		\end{equation*}
		It follows that there exists $ C > 0 $ such that
		\begin{multline*}
			| f((1-u') \langle \rho, \phi \rangle + u \phi(x)) - f((1-u) \langle \rho, \phi \rangle + u \phi(x)) \\ - (u'-u) (\phi(x) - \langle \rho, \phi \rangle) f'(\langle \rho, \phi \rangle) | \leq C (u + u') | u'-u |.
		\end{multline*}
		We then note that the integral of the factor of $ (u'-u) $ on the left hand side with respect to $ \rho(dx) $ is equal to zero, and we replace $ u' $ and $ u $ by $ \frac{y_K}{n + y_K} $ and $ \frac{y}{n+y} $, respectively, to obtain
		\begin{equation} \label{bound_Theta_yK_y}
			\abs{ \Theta_{n, y_K} F_{f,\phi}(\rho) - \Theta_{n,y} F_{f,\phi}(\rho) } \leq \frac{2C}{n K} y.
		\end{equation}
		Integrating with respect to $ y $, we obtain
		\begin{align*}
			b p_0 n \int_{\frac{1}{K}}^{1} \abs{ \Theta_{n, \lfloor K y \rfloor / K} F_{f,\phi}(\rho) - \Theta_{n,y} F_{f,\phi}(\rho) } \frac{dy}{y^{1+\alpha}} &\leq 2 C b p_0 \frac{1}{K} \int_{\frac{1}{K}}^{1} \frac{dy}{y^\alpha} \\
			&= 2 C b p_0 \frac{1}{K} \frac{K^{\alpha-1}-1}{\alpha-1}.
		\end{align*}
		We then see that the right hand side vanishes as $ K \to \infty $ since $ \alpha < 2 $.
		Combined with \eqref{bound_R3_large_y}, this implies the existence of a constant $ C > 0 $ such that
		\begin{equation} \label{bound_R3}
			\abs{ R^{f,\phi}_{K,3}(\rho, n) } \leq \frac{C}{K^{2-\alpha}}.
		\end{equation}
		
		\paragraph*{Bound on $ R^{f,\phi}_{K,4}(\rho, n) $}
		
		Finally, we use \eqref{bound_Theta_small_y} again to write
		\begin{align*}
			\abs{ R^{f,\phi}_{K,4}(\rho, n) } &\leq (d + cn) n K^\alpha \| f'' \| \| \phi \|^2 \left( \frac{\frac{1}{K}}{n - \frac{1}{K}} \right)^2 \\
			&\leq (d + cn) \frac{n}{(n - \frac{1}{K})^2} \| f'' \| \| \phi \|^2 \frac{1}{K^{2-\alpha}}.
		\end{align*}
		It follows that there exists a constant $ C > 0 $ such that, for all $ K $ large enough and $ n \geq c_0 / 2 $,
		\begin{equation} \label{bound_R4}
			\abs{ R^{f,\phi}_{K,4}(\rho, n) } \leq \frac{C}{K^{2-\alpha}}.
		\end{equation}
		
		\paragraph*{Conclusion of the proof}
		
		Combining \eqref{bound_R_K1}, \eqref{bound_R2}, \eqref{bound_R3}, \eqref{bound_R4} and recalling \eqref{diff_AK_A}, we obtain that, for all $ n \geq n_* / 2 $,
		\begin{equation*}
			\sum_{i=1}^{4} \abs{R^{f,\phi}_{K,i}(\rho, n)} \leq \eta_K
		\end{equation*}
		for some sequence $ (\eta_K, K \geq 1) $ which converges to zero as $ K \to \infty $.
		Then, setting
		\begin{equation*}
			\varepsilon^{f,\phi}_K(t) := \int_{0}^{t\wedge \tau_K} \sum_{i=1}^{4} R^{f,\phi}_{K,i}(\rho^K_s, \Nbar(s)) ds,
		\end{equation*}
		we obtain the result, namely that
		\begin{equation*}
			F_{f,\phi}(\rho^K_{t\wedge \tau_K}) - F_{f,\phi}(\rho^K_0) - \int_{0}^{t \wedge \tau_K} A F_{f,\phi}(\rho^K_s, \Nbar(s)) ds - \varepsilon^{f,\phi}_K(t)
		\end{equation*}
		is an $ \F^K_t $-martingale and that
		\begin{equation*}
			\sup_{t \in [0, T]} \abs{ \varepsilon^{f,\phi}_K(t) } \leq \eta_K T
		\end{equation*}
		almost surely.
	\end{proof}
	
	\subsection[Martingale problem for N(t)]{Martingale problem for $ (\Nbar(t), t \geq 0) $} \label{subsec:Nbar}
	
	We now prove Lemma~\ref{lemma:Nbar}, which will conclude the proof of Theorem~\ref{thm:alpha-stable}.
	
	\begin{proof}[Proof of Lemma~\ref{lemma:Nbar}]
		Recall that the infinitesimal generator of $ (\Nbar(t), t \geq 0) $ is given by
		\begin{equation*}
			\mathcal{B}_K g(n) = K^\alpha bn \sum_{k=1}^{\infty} p_k \left( g\left( n + \frac{k}{K} \right) - g(n) \right) + K^\alpha (d + cn) n \left( g\left( n - \frac{1}{K} \right) - g(n) \right).
		\end{equation*}
		Recalling the definition of the operator $ \mathcal{B} $ in \eqref{def:B}, we then write, for $ g = g_\lambda $
		\begin{equation*}
			\mathcal{B}_K g_\lambda(n) = K^{\alpha-1} \mathcal{B} g_\lambda(n) + E^{(1)}_{K,\lambda}(n) + E^{(2)}_{K,\lambda}(n) + E^{(3)}_{K,\lambda}(n),
		\end{equation*}
		where
		\begin{align*}
			E^{(1)}_{K,\lambda}(n) &:= K^\alpha b n \sum_{k=1}^{\infty} p_k \left( g_\lambda\left( n + \frac{k}{K} \right) - g_\lambda(n) - \frac{k}{K} \1{k \leq K} g_\lambda'(n) \right), \\
			E^{(2)}_{K,\lambda}(n) &:= K^\alpha (d + cn) n \left( g_\lambda\left( n - \frac{1}{K} \right) - g_\lambda(n) + \frac{1}{K} g_\lambda'(n) \right), \\
			E^{(3)}_{K,\lambda}(n) &:= - K^{\alpha-1} b n \sum_{k > K} k p_k g_\lambda'(n).
		\end{align*}
		We then bound each term separately.
		
		\paragraph*{Bound on $ E^{(1)}_{K,\lambda} $}
		
		Using \eqref{integral_representation_sum_pk}, we can write
		\begin{equation*}
			E^{(1)}_{K,\lambda}(n) = b n \int_{\frac{1}{K}}^{\infty} \left( g_\lambda \left( n + y_K \right) - g_\lambda (n) - y_K \1{y_K \leq 1} g_\lambda'(n) \right) \frac{p(K y)}{y^{1+\alpha}} dy
		\end{equation*}
		Then, for $ 0 \leq y_K \leq 1 $, by Taylor's formula,
		\begin{equation*}
			g_\lambda \left( n + y_K \right) - g_\lambda (n) - y_K g_\lambda'(n) = (y_K)^2 \int_{0}^{1} (1-u) g_\lambda''(n + u y_K) du.
		\end{equation*}
		By the definition of $ g_\lambda $,
		\begin{align*}
			g_\lambda''(n + u y_K) &= \lambda^2 e^{-\lambda (n + u y_K)} \\
			&\leq \lambda^2 e^{-\lambda n}.
		\end{align*}
		We thus obtain (also using $ y_k \leq y $),
		\begin{equation*}
			0 \leq (y_K)^2 \int_{0}^{1} (1-u) g_\lambda''(n + u y_K) du  \leq y^2 \frac{\lambda^2}{2} e^{-\lambda n}.
		\end{equation*}
		In addition, for $ y_K > 1 $,
		\begin{align*}
			0 \geq g_\lambda \left( n + y_K \right) - g_\lambda (n) &= e^{-\lambda n} \left( e^{-\lambda y_K} - 1 \right) \\
			&\geq - e^{-\lambda n}.
		\end{align*}
		We thus obtain that
		\begin{align*}
			\abs{ E^{(1)}_{K,\lambda}(n) } &\leq C_p b n e^{-\lambda n} \left( \frac{\lambda^2}{2} \int_{\frac{1}{K}}^{1} y^{1-\alpha} dy + \int_{1}^{\infty} \frac{dy}{y^{1+\alpha}} \right) \\
			&\leq C_p b n e^{-\lambda n} \left( \frac{\lambda^2}{2(2-\alpha)} + \frac{1}{\alpha} \right),
		\end{align*}
		where $ C_p := \sup_{y \in [1,\infty)} p(y) $ was defined in \eqref{def:Cp}.
		As a result, there exists a constant $ \newCst{E1} > 0 $ (which depends on $ \lambda $) such that, for all $ n \geq 0 $,
		\begin{equation} \label{bound_E1}
			\abs{E^{(1)}_{K,\lambda}(n)} \leq \Cst{E1}.
		\end{equation}
		
		\paragraph*{Bound on $ E^{(2)}_{K,\lambda} $}
		
		By Taylor's formula,
		\begin{equation*}
			g_\lambda\left( n - \frac{1}{K} \right) - g_\lambda(n) + \frac{1}{K} g_\lambda'(n) = \frac{1}{K^2} \int_{0}^{1} (1-u) g_\lambda''\left( n - \frac{u}{K} \right) du.
		\end{equation*}
		By the definition of $ g_\lambda $,
		\begin{equation*}
			0 \leq g_\lambda''\left( n - \frac{u}{K} \right) \leq \lambda^2 e^{-\lambda n} e^{\lambda / K}.
		\end{equation*}
		As a result,
		\begin{equation*}
			\abs{E^{(2)}_{K,\lambda}(n)} \leq \frac{1}{K^{2-\alpha}} (d + cn) n e^{-\lambda n} \frac{\lambda^2}{2} e^{\lambda / K}.
		\end{equation*}
		Hence there exists a constant $ \newCst{E2} > 0 $ (depending on $ \lambda $) such that, for all $ n \geq 0 $ and $ K \geq 1 $,
		\begin{equation} \label{bound_E2}
			\abs{E^{(2)}_{K,\lambda}(n)} \leq \frac{\Cst{E2}}{K^{2-\alpha}}.
		\end{equation}
		
		\paragraph*{Bound on $ E^{(3)}_{K,\lambda} $}
		
		By \eqref{integral_representation_sum_pk},
		\begin{align*}
			K^{\alpha-1} \sum_{k > K} k p_k &\leq C_p \int_{1}^{\infty} \frac{dy}{y^\alpha} \\
			&= \frac{C_p}{\alpha-1}.
		\end{align*}
		Hence
		\begin{equation*}
			0 \leq E^{(3)}_{K,\lambda}(n) \leq \lambda \frac{b C_p}{\alpha-1} n e^{-\lambda n}.
		\end{equation*}
		Thus there exists a constant $ \newCst{E3} > 0 $ (depending on $ \lambda $) such that, for all $ n \geq 0 $,
		\begin{equation} \label{bound_E3}
			\abs{E^{(3)}_{K,\lambda}(n)} \leq \Cst{E3}.
		\end{equation}
		
		\paragraph*{Conclusion of the proof}
		
		Setting
		\begin{equation*}
			\delta^\lambda_K(t) := - \int_{0}^{t} \sum_{i=1}^{3} E^{(i)}_{K,\lambda}(\Nbar(s)) ds,
		\end{equation*}
		we obtain that
		\begin{equation*}
			g_\lambda(\Nbar(t)) - g_\lambda(\Nbar(0)) - K^{\alpha-1} \int_{0}^{t} \mathcal{B} g_\lambda(\Nbar(s)) ds + \delta^\lambda_K(t)
		\end{equation*}
		is an $ \mathcal{F}^K_t $-martingale, and
		\begin{equation*}
			\sum_{t \in [0,T]} \abs{ \delta^\lambda_K(t) } \leq \left( \Cst{E1} + \frac{\Cst{E2}}{K^{2-\alpha}} + \Cst{E3} \right) T,
		\end{equation*}
		which concludes the proof of the lemma.
	\end{proof}
	
	\section[Convergence to a Lambda-Fleming-Viot process for Neveu's logistic branching process]{Convergence to a $ \Lambda $-Fleming-Viot process for Neveu's logistic branching process} \label{sec:neveu}
	
	In this section, we prove Theorem~\ref{thm:neveu}.
	The overall strategy is similar to that used in the previous section to prove Theorem~\ref{thm:alpha-stable}, although some additional care is needed due to the fact that the offspring distribution does not have a finite first moment and because the scaling of the population size is different.
	
	\subsection{Proof of the main result}
	
	Recall that, in the setting of Theorem~\ref{thm:neveu},
	\begin{align*}
		\Nbar(t) := \frac{1}{K \log(K)} \langle \nu^K_t, 1 \rangle, && \rho^K_t := \varrho(\nu^K_t).
	\end{align*}
	Let $ (\F^K_t, t \geq 0) $ denote the natural filtration associated to $ \lbrace (\rho^K_t, \Nbar(t)), t \geq 0 \rbrace $.
	By the same argument as the one given in Section~\ref{subsec:proof_stable}, we can assume, without loss of generality, that there exist $ c_0 \in (0, n_*) $ and $ C_0 > 0 $ such that, for all $ K \geq 1 $,
	\begin{equation} \label{N0_compact_neveu}
		c_0 \leq \Nbar(0) \leq C_0,
	\end{equation}
	almost surely.
	We then define a stopping time $ \tau_K $ as
	\begin{equation*}
		\tau_K := \inf \left\lbrace t \geq 0 : \Nbar(t) < \frac{c_0}{2} \right\rbrace.
	\end{equation*}
	In Subsection~\ref{subsec:tau_K_neveu}, we prove the following.
	
	\begin{lemma} \label{lemma:tauK_neveu}
		Under the assumptions of Theorem~\ref{thm:neveu}, for any $ t \geq 0 $,
		\begin{equation*}
			\lim_{K \to \infty} \P{ \tau_K \leq t } = 0.
		\end{equation*}
	\end{lemma}
	
	In Subsection~\ref{subsec:Nbar_neveu}, we prove the following, which allows us to conclude that the occupation measure of $ (\Nbar(t), t \geq 0) $ concentrates on $ \lbrace n_* \rbrace $ as $ K \to \infty $.
	Let $ \mathcal{B} : C^1(\R_+) \to C(\R_+) $ be the operator defined as
	\begin{equation} \label{def:B_neveu}
		\mathcal{B} g(n) := c(n_* - n) n g'(n).
	\end{equation}
	As in Section~\ref{sec:alpha_stable}, for $ \lambda > 0 $, let $ g_\lambda : \R_+ \to [0,1] $ be the function defined by
	\begin{equation*}
		g_\lambda(n) := e^{-\lambda n}.
	\end{equation*}
	
	\begin{lemma} \label{lemma:Nbar_neveu}
		Under the assumptions of Theorem~\ref{thm:neveu} and \eqref{N0_compact_neveu}, for any $ \lambda > 0 $, there exists a process $ (\delta^\lambda_K(t), t \geq 0) $ such that
		\begin{equation*}
			g_\lambda(\Nbar(t)) - g_\lambda(\Nbar(0)) - \log(K) \int_{0}^{t} \mathcal{B} g_\lambda(\Nbar(s)) ds + \delta^\lambda_K(t)
		\end{equation*}
		is an $ \F^K_t $-martingale.
		In addition, for any $ T > 0 $, there exists a deterministic sequence $ (\eta_K, K \geq 1) $ such that, for all $ K \geq 1 $,
		\begin{equation*}
			\sup_{t \in [0,T]} \abs{\delta^\lambda_K(t)} \leq \eta_K,
		\end{equation*}
		almost surely, and $ \eta_K = \littleO{\log(K)} $.
	\end{lemma}
	
	We also need a bound similar to that in \eqref{lemma:bounded_expectation} to show that the family of random variables $ \lbrace \Nbar(t \wedge \tau_K), t \in [0,T], K \geq 1 \rbrace $ is tight.
	Here, however, $ \Nbar(t) $ does not have a finite expectation, so we need a different estimate. 
	This turns out to be quite delicate, as we shall see below.
	For $ \varepsilon \in (0, n_*) $, let $ V_\varepsilon : (0,\infty) \to \R $ be defined as
	\begin{equation*}
		V_\varepsilon(n) := \frac{n}{n_* + \varepsilon} - 1 - \log\left( \frac{n}{n_* - \varepsilon} \right).
	\end{equation*}
	where $ n_* $ is given by \eqref{def:nbar}.
	We then have the following result, which is proved in Subsection~\ref{subsec:proof_moment_bound_neveu}.
	
	\begin{lemma} \label{lemma:bound_E_log_V_n}
		Under the assumptions of Theorem~\ref{thm:neveu}, for any $ T > 0 $, and $ \varepsilon > 0 $ such that
		\begin{equation} \label{lower_bound_Vepsilon}
			\inf_{n \geq 0} V_\varepsilon(n) > -1,
		\end{equation}
		there exists $ C_{T,\varepsilon} > 0 $ such that
		\begin{equation*}
			\sup_{K \geq 1} \sup_{t \in [0,T]} \E{\log(1+ V_\varepsilon(\Nbar(t \wedge \tau_K)))} \leq C_{T,\varepsilon}.
		\end{equation*}
	\end{lemma}
	
	Note that \eqref{lower_bound_Vepsilon} is satisfied for $ \varepsilon > 0 $ small enough since
	\begin{equation*}
		\inf_{n \geq 0} V_\varepsilon(n) = V_\varepsilon(n_* + \varepsilon) = - \log\left( 1 + \frac{2 \varepsilon}{n_* - \varepsilon} \right).
 	\end{equation*}
	
	To complete the proof of Theorem~\ref{thm:neveu}, we need the following lemma, which allows us to characterise the possible limits of $ (\rho^K_t, t \geq 0) $.
	Let $ A : \mathcal{D}_0 \subset C(\mathcal{M}_1(\X)) \to C(\mathcal{M}_1(\X)) $ be the linear operator defined by
	\begin{equation} \label{def:A_neveu}
		A F_{f, \phi} (\rho) := b p_0 \int_{0}^{1} \int_\X \left( f\left( (1-u) \langle \rho, \phi \rangle + u \phi(x) \right) - f(\langle \rho, \phi \rangle) \right) \rho(dx) du,
	\end{equation}
	and extend $ A $ to $ \mathcal{D} $ by linearity.
	We recognise here the generator of the $ \Lambda $-Fleming-Viot process with $ \Lambda(du) = b p_0 du $.
	
	\begin{lemma} \label{lemma:neutral_types_neveu}
		Under the assumptions of Theorem~\ref{thm:neveu} and \eqref{N0_compact_neveu}, for any $ F = F_{f,\phi} \in \mathcal{D}_0 $, there exists a process $ (\varepsilon^{f,\phi}_K(t), t \geq 0) $ such that
		\begin{equation*}
			F_{f,\phi}(\rho^K_{t \wedge \tau_K}) - F_{f,\phi}(\rho^K_0) - \int_{0}^{t \wedge \tau_K} A F_{f,\phi}(\rho^K_s) ds + \varepsilon^{f,\phi}_K(t)
		\end{equation*}
		is an $ \F^K_t $-martingale.
		In addition, for any $ T \geq 0 $, there exists a deterministic sequence $ (\eta_K, K \geq 1) $, converging to zero as $ K \to \infty $ and such that, for all $ K \geq 1 $,
		\begin{equation*}
			\sup_{t \in [0, T]} \abs{ \varepsilon^{f,\phi}_K(t) } \leq \eta_K, \quad \text{ almost surely.}
		\end{equation*}
	\end{lemma}
	
	We prove Lemma~\ref{lemma:neutral_types_neveu} in Subsection~\ref{subsec:proof_types_neveu}.
	We can now prove Theorem~\ref{thm:neveu}.
	
	\begin{proof}[Proof of Theorem~\ref{thm:neveu}]
		The proof proceeds along the same steps as that of Theorem~\ref{thm:alpha-stable}, applying Theorem~\ref{thm:averaging} and Proposition~\ref{prop:averaging_Y} to $ X_K = \rho^K $ and $ Y_K = \Nbar $.
		Let $ \Gamma_K $ denote the random measure taking values in $ \ell(\R_+) $ defined as
		\begin{equation*}
			\Gamma_K([0,t] \times B) := \int_{0}^{t \wedge \tau_K} \1{\Nbar(s) \in B} ds,
		\end{equation*}
		for all measurable $ B \subset \R_+ $ and $ t \geq 0 $.
		
		Assumption~\ref{ass:compact_containment} is again satisfied since $ \M $ is compact.
		Assumption~\ref{ass:martingale_pb} follows by Lemma~\ref{lemma:neutral_types_neveu} by linearity, noting that, for any $ F \in \mathcal{D}_0 $, there exists a constant $ C > 0 $ such that, for all $ \rho \in \M $,
		\begin{equation*}
			| AF(\rho) | \leq C.
		\end{equation*}
		Assumption~\ref{ass:tightness_Y} is satisfied by Lemma~\ref{lemma:bound_E_log_V_n}.
		Finally, $ \mathcal{D} $ is dense in $ C(\M) $ and Assumption~\ref{ass:stopping_time} follows by Lemma~\ref{lemma:tauK_neveu}.
		
		As a result, by Theorem~\ref{thm:averaging}, $ \lbrace (\rho^K, \Gamma_K), K \geq 1 \rbrace $ is tight in $ D([0,\infty), \M) \times \ell(\R_+) $ and, for any limit point $ (\rho, \Gamma) $, for all $ F \in \mathcal{D} $,
		\begin{equation*}
			F(\rho_t) - F(\rho_0) - \int_{0}^{t} AF(\rho_s) ds
		\end{equation*}
		is an $ \F_t $-martingale, where
		\begin{equation*}
			\F_t := \sigma \lbrace \rho_s, \Gamma([0,s] \times B), s \in [0,t], B \in \mathcal{B}(\R_+) \rbrace.
		\end{equation*}
		This shows that any limit point $ \rho $ of the sequence $ (\rho^K, K\geq 1) $ is a solution to the martingale problem associated to $ (A, \mathcal{D}) $, and yields the convergence in distribution of $ \rho^K $ as $ K \to \infty $.
			
		In addition, by Lemma~\ref{lemma:Nbar_neveu}, we can apply Proposition~\ref{prop:averaging_Y} with $ \beta_K = \log(K) $ to obtain that, for any $ \lambda > 0 $,
		\begin{equation*}
			\int_{[0,t] \times \R_+} (b p_0 - c n) n e^{-\lambda n} \Gamma(ds, dn) = 0, \quad \forall t \geq 0,
		\end{equation*}
		almost surely.
		It follows that $ \Gamma([0,t], dn) $ is concentrated on $ \lbrace 0, n_* \rbrace $.
		In addition, by the definition of $ \Gamma_K $ and $ \tau_K $, for any continuous $ f : \R_+ \to \R $ supported on $ [0, c_0 / 2) $,
		\begin{equation*}
			\int_{[0,t] \times \R_+} f(n) \Gamma(ds, dn) = 0,
		\end{equation*}
		almost surely.
		Combined with \eqref{Gamma_t}, this shows that $ \Gamma(ds, dn) = \delta_{n_*}(dn) ds $, almost surely, which entails \eqref{concentration_Nbar_neveu}.
		This concludes the proof of Theorem~\ref{thm:neveu}.
	\end{proof}
	
	\subsection{Control of the stopping time} \label{subsec:tau_K_neveu}
	
	Let us now prove Lemma~\ref{lemma:tauK_neveu}.
	We start by a technical lemma.
	
	\begin{lemma} \label{lemma:cvg_n_K}
		Let $ (b_K, K \geq 1) $ be a sequence of integers such that
		\begin{equation} \label{log_equivalent_bK}
			\lim_{K \to \infty} \frac{\log(b_K)}{\log(K)} = \gamma \in (0,\infty),
		\end{equation}
		then
		\begin{equation*}
			\frac{1}{c \log(K)} \left( b \sum_{k=1}^{b_K} k p_k - d \right) \cvgas{K} \gamma n_*.
		\end{equation*}
	\end{lemma}
	
	\begin{proof}
		Set
		\begin{equation*}
			n_K := \frac{1}{c \log(K)} \left( b \sum_{k=1}^{b_K} k p_k - d \right).
		\end{equation*}
		By the definition of $ n_* $,
		\begin{equation*}
			n_K - \gamma n_* = \frac{b}{c \log(K)} \left( \sum_{k=1}^{b_K} k p_k - \gamma p_0 \log(K) \right) - \frac{d}{c \log(K)}.
		\end{equation*}
		We then write
		\begin{equation} \label{log_equiv}
			\sum_{k=1}^{b_K} k p_k - \gamma p_0 \log(K) = \sum_{k=1}^{b_K} (k^2 p_k - p_0) \frac{1}{k} + p_0 \left( \sum_{k=1}^{b_K} \frac{1}{k} - \gamma \log(K) \right).
		\end{equation}
		Since
		\begin{equation*}
			\sum_{k=1}^{b_K} \frac{1}{k} = \log(b_K) + \bigO{1},
		\end{equation*}
		the second term on the right of \eqref{log_equiv} is
		\begin{align*}
			\log(b_K) - \gamma \log(K) + \bigO{1} &= \log(K) \left( \frac{\log(b_K)}{\log(K)} - \gamma \right) + \bigO{1} \\
			&= \littleO{\log(K)},
		\end{align*}
		by \eqref{log_equivalent_bK}.
		We now bound the first term on the right of \eqref{log_equiv}.
		To do so, take $ \varepsilon > 0 $ and let $ k_\varepsilon \in \N $ be such that, for all $ k \geq k_\varepsilon $, $ | k^2 p_k - p_0 | \leq \varepsilon $.
		Then, for $ b_K \geq k_\varepsilon $
		\begin{equation*}
		 \left| \sum_{k=1}^{b_K} (k^2 p_k - p_0) \frac{1}{k} \right| \leq C \sum_{k=1}^{k_\varepsilon-1} \frac{1}{k} + \varepsilon \sum_{k=k_\varepsilon}^{b_K} \frac{1}{k} = \varepsilon \log(b_K) + \bigO{1}.
		\end{equation*}
		As a result, using \eqref{log_equivalent_bK},
		\begin{equation*}
			\limsup_{K \to \infty} | n_K - \gamma n_* | \leq \gamma \varepsilon,
		\end{equation*}
		for all $ \varepsilon > 0 $.
		Letting $ \varepsilon \downarrow 0 $, we obtain the result.
	\end{proof}

	Let $ \mathcal{B}_K $ denote the infinitesimal generator of $ (\Nbar(t), t \geq 0) $, i.e., for $ g : \R_+ \to \R $ bounded and measurable, setting
	\begin{equation*}
	a_K := K \log(K),
	\end{equation*}
	we have
	\begin{multline} \label{def:BK_neveu}
		\mathcal{B}_K g(n) = a_K b n \sum_{k=1}^{+\infty} p_k \left( g\left( n + \frac{k}{a_K} \right) - g(n) \right)  \\ + a_K \left( d + c \log(K) n \right) n \left( g\left( n - \frac{1}{a_K} \right) - g(n) \right).
	\end{multline}
	We note that, by a slight adaptation of Proposition~\ref{prop:coupling}, for any $ \beta \in (0,1) $ such that
	\begin{equation*}
		\beta n_* \in \left( \frac{c_0}{2}, c_0 \right),
	\end{equation*}
	we can construct a process $ (\Nbeta(t), t \geq 0) $ taking values in $ \frac{1}{a_K} \N $ such that
	\begin{itemize}
		\item $ \Nbeta(0) $ converges to $ \beta n_* $ in probability as $ K \to \infty $,
		
		\item $ (\Nbeta(t), t \geq 0) $ is a Markov process with generator $ \widetilde{\mathcal{B}}_K $ defined as
		\begin{multline*}
			\widetilde{\mathcal{B}}_K g(n) := a_K b n \sum_{k=1}^{\lfloor K^\beta \rfloor} p_k \left( g\left( n + \frac{k}{a_K} \right) - g(n) \right) \\ + a_K (d + c \log(K) n) n \left( g\left( n - \frac{1}{a_K} \right) - g(n) \right),
		\end{multline*}
		
		\item and $ \Nbeta(t) \leq \Nbar(t) $ for all $ t \geq 0 $, almost surely.
	\end{itemize}
	Let us then define $ \tilde{n}_K $ as
	\begin{equation*}
		\tilde{n}_K := \frac{1}{c \log(K)} \left( b \sum_{k=1}^{\lfloor K^\beta \rfloor} k p_k - d \right).
	\end{equation*}
	Lemma~\ref{lemma:cvg_n_K} implies that
	\begin{equation} \label{cvg_n_tilde}
		\lim_{K \to \infty} \tilde{n}_K = \beta n_*.
	\end{equation}
	For any $ \varepsilon \in (0, \beta n_*) $, let us then define a stopping time $ \tilde{\tau}_K $ as
	\begin{equation*}
		\tilde{\tau}_K := \inf \lbrace t \geq 0 : | \Nbeta(t) - \tilde{n}_K | > \varepsilon \rbrace.
	\end{equation*}
	
	\begin{lemma} \label{lemma:tau_epsilon}
		Under the assumptions of Theorem~\ref{thm:neveu} and \eqref{N0_compact_neveu}, for any $ \varepsilon \in (0, \beta n_*) $ and any $ t \geq 0 $,
		\begin{equation*}
			\lim_{K \to \infty} \P{\tilde{\tau}_K \leq t} = 0.
		\end{equation*}
	\end{lemma}
	
	Lemma~\ref{lemma:tauK_neveu} then follows directly from \eqref{cvg_n_tilde} and Lemma~\ref{lemma:tau_epsilon} by the following considerations.
	
	\begin{proof}[Proof of Lemma~\ref{lemma:tauK_neveu}]
		Since $ \beta n_* \in \left( \frac{c_0}{2}, c_0 \right) $, we can find $ \varepsilon > 0 $ such that $ \beta n_* - 2\varepsilon > c_0/2 $.
		By \eqref{cvg_n_tilde}, for all $ K $ large enough, $ \tilde{n}_K \geq \beta n_* - \varepsilon $.
		Then, by the definition of $ \tau_K $ and $ \tilde{\tau}_K $ and the fact that $ \Nbeta(t) \leq \Nbar(t) $ for all $ t \geq 0 $, $ \tilde{\tau}_K \leq \tau_K $ almost surely, and Lemma~\ref{lemma:tauK_neveu} follows by Lemma~\ref{lemma:tau_epsilon}.
	\end{proof}
	
	We now conclude this section by proving Lemma~\ref{lemma:tau_epsilon}.
	
	\begin{proof}[Proof of Lemma~\ref{lemma:tau_epsilon}]
		Let $ g_K : \R_+ \to \R_+ $ be defined as
		\begin{equation*}
			g_K(n) := \frac{n}{\tilde{n}_K} - 1 - \log\left( \frac{n}{\tilde{n}_K} \right).
		\end{equation*}
		Recall from \eqref{quadratic_bound_V} that, for all $ n \neq \tilde{n}_K $, $ g_K(n) > g_K(\tilde{n}_K) = 0 $ and, by \eqref{cvg_n_tilde}, there exists $ C_1 > 0 $ and $ C_2 > 0 $ such that, for all $ K $ large enough and all $ n \in [\tilde{n}_K - \varepsilon, \tilde{n}_K + \varepsilon] $,
		\begin{equation} \label{quadratic_bound_gK}
			C_1 (n - \tilde{n}_K)^2 \leq g_K(n) \leq C_2 (n - \tilde{n}_K)^2.
		\end{equation}
		We then prove that, for any $ t \geq 0 $,
		\begin{equation} \label{bound_gK}
			\lim_{K \to \infty} \E{ g_K(\Nbeta(t\wedge \tilde{\tau}_K)) } = 0.
		\end{equation}
		To do so, we note that $ \widetilde{\mathcal{B}}_K g_K $ can be written as
		\begin{equation*}
			\widetilde{\mathcal{B}}_K g_K(n) = c \log(K) (\tilde{n}_K - n) n g_K'(n) + \widetilde{B}_{1,K}(n) + \widetilde{B}_{2,K}(n),
		\end{equation*}
		where
		\begin{align*}
			\widetilde{B}_{1,K}(n) &= \frac{1}{a_K} b n \sum_{k=1}^{\lfloor K^\beta \rfloor} k^2 p_k \int_{0}^{1} (1-u) g_K''\left( n + u \frac{k}{a_K} \right) du \\
			\widetilde{B}_{2,K}(n) &= \frac{1}{a_K} (d + c \log(K) n) n \int_{0}^{1} (1-u) g_K''\left( n - \frac{u}{a_K} \right) du.
		\end{align*}
		Then we note that, by the definition of $ g_K $,
		\begin{equation*}
			(\tilde{n}_K - n) n g_K'(n) = - \frac{c}{\tilde{n}_K} (n-\tilde{n}_K)^2 \leq 0.
		\end{equation*}
		In addition, since $ g_K''(n) = \frac{1}{n^2} $, we see that there exists a constant $ C > 0 $ such that, for all $ n \in [\tilde{n}_K - \varepsilon, \tilde{n}_K + \varepsilon] $ and $ k \geq 1 $,
		\begin{equation*}
			\abs{ b n  \int_{0}^{1} (1-u) g_K''\left( n + u \frac{k}{a_K} \right) du } \leq C.
		\end{equation*}
		By \eqref{condition_pk_neveu}, there exists $ C > 0 $ such that $ k^2 p_k \leq C $ for all $ k \in \N $. 
		As a result, there exists $ C > 0 $ such that, for all $ n \in [\tilde{n}_K - \varepsilon, \tilde{n}_K + \varepsilon] $, 
		\begin{equation*}
			\abs{ \widetilde{B}_{1,K}(n) } \leq \frac{C}{K^{1-\beta} \log(K)}.
		\end{equation*}
		By a similar argument, there exists a constant $ C > 0 $ such that, for all $ n \in [\tilde{n}_K - \varepsilon, \tilde{n}_K + \varepsilon] $,
		\begin{equation*}
			\abs{ \widetilde{B}_{2,K}(n) } \leq \frac{C}{K} \left( \frac{d}{\log(K)} + c(\tilde{n}_K+\varepsilon) \right) = \littleO{ \frac{1}{K^{1-\beta} \log(K)} }.
		\end{equation*}
		Hence there exists a constant $ C_3 > 0 $ such that, for all $ n \in [\tilde{n}_K - \varepsilon, \tilde{n}_K + \varepsilon] $,
		\begin{equation} \label{bound_Btilde}
			\widetilde{\mathcal{B}}_K g_K(n) \leq \frac{C_3}{K^{1-\beta} \log(K)}.
		\end{equation}
		We then note that, for all $ K $,
		\begin{equation*}
			Z_K(t) := g_K(\Nbeta(t)) - g_K(\Nbeta(0)) - \int_{0}^{t} \widetilde{\mathcal{B}}_K g_K(\Nbeta(s)) ds
		\end{equation*}
		is a local martingale and, by \eqref{N0_compact},
		\begin{equation*}
			\sup_{K \geq 1} \E{g_K(\Nbeta(0))} \leq \sup_{K \geq 1} \E{g_K(\Nbar(0))} < \infty.
		\end{equation*}
		Since, by \eqref{bound_Btilde}, $ Z_K(t) $ is square-integrable for all $ K $,
		\begin{equation*}
			\E{ Z_K(t\wedge \tilde{\tau}_K) } = 0.
		\end{equation*}
		Using \eqref{bound_Btilde}, this yields
		\begin{equation*}
			\E{ g_K(\Nbeta(t \wedge \tilde{\tau}_K)) } \leq \E{ g_K(\Nbeta(0)) } + \frac{C_3 t}{K^{1-\beta} \log(K)},
		\end{equation*}
		for all $ t \geq 0 $.
		Combining \eqref{cvg_n_tilde} and the fact that $ \Nbeta(0) \to \beta n_* $ in probability, we obtain that
		\begin{equation*}
			\lim_{K \to \infty} \E{ g_K(\Nbeta(0)) } = 0.
		\end{equation*}
		We thus obtain \eqref{bound_gK}.
		To conclude, by \eqref{quadratic_bound_gK},
		\begin{equation*}
			\P{\tilde{\tau}_K \leq t} \leq \frac{\E{ g_K(\Nbeta(t \wedge \tilde{\tau}_K)) }}{C_1 \varepsilon^2}.
		\end{equation*}
		The right hand side vanishes as $ K \to \infty $ by \eqref{bound_gK}, which concludes the proof of Lemma~\ref{lemma:tau_epsilon}.
	\end{proof}
	
	\subsection{The rescaled population process} \label{subsec:Nbar_neveu}
	
	In this subsection, we prove Lemma~\ref{lemma:Nbar_neveu}.
	Recall that the infinitesimal generator of $ (\Nbar(t), t \geq 0) $ is given by
	\begin{multline*}
		\mathcal{B}_K g(n) = a_K b n \sum_{k=1}^{+\infty} p_k \left( g\left( n + \frac{k}{a_K} \right) - g(n) \right)  \\ + a_K \left( d + c \log(K) n \right) n \left( g\left( n - \frac{1}{a_K} \right) - g(n) \right).
	\end{multline*}
	As in Subsection~\ref{subsec:Nbar}, for $ g $ bounded and twice continuously differentiable, we split the sum over $ k $ in two parts and write a Taylor expansion of $ g $ in the first part of the sum to write
	\begin{equation} \label{decomp:BK_neveu}
		\mathcal{B}_K g(n) = \log(K) \mathcal{B} g(n) + \sum_{i=1}^{4} B^K_i g(n),
	\end{equation}
	where $ \mathcal{B}	g(n) = c(n_* - n) n g'(n) $ was defined in \eqref{def:B_neveu} and
	\begin{align*}
		B^K_1 g(n) &:=  c \log(K) \left( \frac{b}{c \log(K)} \sum_{k=1}^{\lfloor a_K \rfloor} k p_k - \frac{d}{c \log(K)} - n_* \right) n g'(n), \\
		B^K_2 g(n) &:= a_K b n \sum_{k > a_K} p_k \left( g\left( n + \frac{k}{a_K} \right) - g(n) \right), \\
		B^K_3 g(n) &:= a_K b n \sum_{k=1}^{\lfloor a_K \rfloor} p_k \left( g\left( n + \frac{k}{a_K} \right) - g(n) - \frac{k}{a_K} g'(n) \right), \\
		B^K_4 g(n) &:= a_K \left( d + c \log(K) n \right) n \left( g\left( n - \frac{1}{a_K} \right) - g(n) + \frac{1}{a_K} g'(n) \right),
	\end{align*}
	(note that the cut-off in the sum is different than that in Subsection~\ref{subsec:Nbar}).
	
	Let us set
	\begin{equation} \label{def:n_bar_K}
		\overline{n}_K := \frac{b}{c \log(K)} \sum_{k=1}^{\lfloor a_K \rfloor} k p_k - \frac{d}{c \log(K)}.
	\end{equation}
	By Lemma~\ref{lemma:cvg_n_K},
	\begin{equation} \label{cvg_n_bar_K}
		\lim_{K \to\infty} \overline{n}_K = n_*.
	\end{equation}
	We now prove Lemma~\ref{lemma:Nbar_neveu}.
	
	\begin{proof}[Proof of Lemma~\ref{lemma:Nbar_neveu}]
		First, for any $ \lambda > 0 $,
		\begin{equation*}
			B^K_1 g_{\lambda}(n) = - c \log(K) (\overline{n}_K - n_*) \lambda n e^{- \lambda n}.
		\end{equation*}
		Since $ n \mapsto \lambda n e^{-\lambda n} $ is bounded, by \eqref{cvg_n_bar_K}, it follows that
		\begin{equation} \label{bound_B1_glambda}
			\sup_{n \geq 0} | B^K_1 g_\lambda (n) | = \littleO{\log(K)}.
		\end{equation}
		Turning to $ B^K_2 $, we note that
		\begin{equation*}
			B^K_2 g_\lambda(n) = b n e^{-\lambda n} a_K \sum_{k > a_K} \left( e^{-\lambda k / a_K} - 1 \right) p_k.
		\end{equation*}
		Hence, for any $ \lambda > 0 $, there exists a constant $ C > 0 $ such that, for all $ n \geq 0 $,
		\begin{equation*}
			| B^K_2 g_\lambda(n) | \leq C a_K \sum_{k > a_K} p_k.
		\end{equation*}
		By \eqref{condition_pk_neveu}, there exists a constant $ C > 0 $ such that
		\begin{align}
			a_K \sum_{k > a_K} p_k \leq C a_K \sum_{k > a_K} \frac{1}{k^2}.
		\end{align}
		Since the right hand side is $ \bigO{1} $, we obtain that, for some constant $ C > 0 $,
		\begin{equation} \label{bound_B2_glambda}
			\sup_{n \geq 0} | B^K_2 g_\lambda (n) | \leq C.
		\end{equation}
		Turning to $ B^K_3 $, we write
		\begin{equation*}
			B^K_3 g_\lambda(n) = a_K b n e^{-\lambda n} \sum_{k=1}^{\lfloor a_K \rfloor} p_k \left( e^{-\lambda k / a_K} - 1 + \frac{\lambda k}{a_K}  \right).
		\end{equation*}
		By Taylor's formula,
		\begin{equation*}
			e^{-\lambda k / a_K} - 1 + \frac{\lambda k}{a_K} = \left( \frac{\lambda k}{a_K} \right)^2 \int_{0}^{1} (1-u) e^{-u \lambda k / a_K } du.
		\end{equation*}
		As a result, for any $ \lambda > 0 $, there exists $ C > 0 $ such that
		\begin{equation*}
			\sup_{n \geq 0} | B^K_3 g_\lambda(n) | \leq \frac{C}{a_K} \sum_{k=1}^{\lfloor a_K \rfloor} k^2 p_k.
		\end{equation*}
		By \eqref{condition_pk_neveu}, the right hand side is bounded by a constant independent of $ K $, hence there exists $ C > 0 $ such that
		\begin{equation} \label{bound_B3_glambda}
			\sup_{n \geq 0} | B^K_3 g_\lambda(n) | \leq C.
		\end{equation}
		Finally, using Taylor's formula again, we write
		\begin{equation*}
				B^K_4 g_\lambda (n) = (d + c \log(K) n) n e^{-\lambda n} \frac{\lambda^2}{a_K} \int_{0}^{1} (1-u) e^{u \lambda / a_K} du.
		\end{equation*}
		Since both $ n \mapsto n e^{-\lambda n} $ and $ n \mapsto n^2 e^{-\lambda n} $ are bounded on $ (0,\infty) $, for any $ \lambda > 0 $, there exists $ C > 0 $ such that, for all $ n \geq 0 $,
		\begin{equation*}
			| B^K_4 g_\lambda(n) | \leq C \left( \frac{d}{a_K} + \frac{c}{K} \right).
		\end{equation*}
		Hence there exists a constant $ C > 0 $ such that
		\begin{equation} \label{bound_B4_glambda}
			\sup_{n \geq 0} | B^K_4 g_\lambda(n) | \leq \frac{C}{K}.
		\end{equation}
		To conclude, we set
		\begin{equation*}
			\delta^\lambda_K(t) = - \int_{0}^{t} \sum_{i=1}^{4} B^K_i g_\lambda(\Nbar(s)) ds,
		\end{equation*}
		and, combining \eqref{bound_B1_glambda}, \eqref{bound_B2_glambda}, \eqref{bound_B3_glambda} and \eqref{bound_B4_glambda}, we obtain that, for any $ T > 0 $ and any $ \lambda > 0 $, there exists a deterministic sequence $ (\eta_K, K > 0) $ such that, almost surely,
		\begin{equation*}
		 	\sup_{t \in [0, T]} | \delta^\lambda_K(t) | \leq \eta_K,
		\end{equation*}
		and $ \eta_K = \littleO{\log(K)} $.
		This concludes the proof of the lemma.
	\end{proof}
	
	\subsection{Stochastic averaging for the distribution of neutral types} \label{subsec:proof_types_neveu}
	
	In this subsection, we prove Lemma~\ref{lemma:neutral_types_neveu}.
	Recall the notation introduced in \eqref{def:Theta}, and define $ A_K : \mathcal{D}_0 \to C(\mathcal{M}_1(\X) \times \R_+) $ as
	\begin{equation*}
		A_K F_{f,\phi}(\rho, n) = a_K b n \sum_{k=1}^{\infty} p_k \Theta_{n, k/a_K} F_{f,\phi}(\rho) + a_K (d + c \log(K) n) n \, \Theta_{n, -1/a_K} F_{f,\phi}(\rho).
	\end{equation*}
	Then, by definition, 
	\begin{equation*}
		F_{f,\phi}(\rho^K_t) - F_{f,\phi}(\rho^K_0) - \int_{0}^{t} A_K F_{f,\phi}(\rho^K_s, \Nbar(s)) ds
	\end{equation*}
	is a local martingale with respect to the filtration $ (\F^K_t, t \geq 0) $.
	We can now prove Lemma~\ref{lemma:neutral_types_neveu}.
	
	\begin{proof}[Proof of Lemma~\ref{lemma:neutral_types_neveu}]
		Using \eqref{def:p_y}, we write
		\begin{align*}
			A_K F_{f,\phi}(\rho, n) &= a_K b n \int_{1}^{\infty} \Theta_{n, \lfloor y \rfloor / a_K} F_{f,\phi}(\rho) \frac{p(y)}{y^2} dy \\ &\hspace{4cm} + a_K (d + c \log(K) n) n \, \Theta_{n, -1/a_K} F_{f,\phi}(\rho) \\
			&= b n \int_{1/a_K}^{\infty} \Theta_{n, y_K} F_{f,\phi}(\rho) \frac{p(a_K y)}{y^2} dy +  a_K (d + c \log(K) n) n \, \Theta_{n, -1/a_K} F_{f,\phi}(\rho),
		\end{align*}
		where $ y_K := \lfloor a_K y \rfloor / a_K $.
		In parallel, we note that, for any $ n > 0 $, by a change of variables,
		\begin{equation*}
			b p_0 n \int_{0}^{\infty} \Theta_{n,y} F_{f,\phi}(\rho) \frac{dy}{y^2} = A F_{f,\phi}(\rho).
		\end{equation*}
		Subtracting the above expressions, we obtain
		\begin{equation*}
			A_K F_{f,\phi}(\rho, n) - A F_{f,\phi}(\rho) = \sum_{i=1}^{4} R^{f,\phi}_{K,i}(\rho, n),
		\end{equation*}
		where
		\begin{align*}
			R^{f,\phi}_{K,1}(\rho, n) &= a_K (d + c \log(K) n) n \,  \Theta_{n, -1/a_K} F_{f,\phi}(\rho), \\
			R^{f,\phi}_{K,2}(\rho, n) &= - b p_0 n \int_{0}^{1/a_K} \Theta_{n,y} F_{f,\phi}(\rho) \frac{dy}{y^2}, \\
			R^{f,\phi}_{K,3}(\rho, n) &= b n \int_{1/a_K}^{\infty}  \Theta_{n, y_K} F_{f,\phi}(\rho)  \frac{ p(a_K y) - p_0 }{y^2} dy, \\
			R^{f,\phi}_{K,4}(\rho, n) &= b p_0 n \int_{1/a_K}^{\infty} \left( \Theta_{n, y_K} F_{f,\phi}(\rho) - \Theta_{n, y} F_{f,\phi}(\rho) \right) \frac{dy}{y^2}.
		\end{align*}
		We then bound each term separately.
		
		\paragraph*{Bound on $ R^{f,\phi}_{K,1} $}
		
		Using \eqref{bound_Theta_small_y},
		\begin{equation*}
			|R^{f,\phi}_{K,1} (\rho, n)| \leq \frac{1}{a_K} \| f'' \| \| \phi \|^2 (d + c \log(K) n) \frac{n}{(n - \frac{1}{a_K})^2}.
		\end{equation*}
		As a result, there exists a constant $ C > 0 $ such that, for all $ n \geq c_0 / 2 $ and $ K $ large enough,
		\begin{equation} \label{bound_RK1_neveu}
			|R^{f,\phi}_{K,1}(\rho, n)| \leq \frac{C}{K}.
		\end{equation}
		
		\paragraph*{Bound on $ R^{f,\phi}_{K,2} $}
		
		Using \eqref{bound_Theta_small_y}, we obtain
		\begin{align*}
			|R^{f,\phi}_{K,2}(\rho, n)| &\leq b p_0 \| f'' \| \| \phi \|^2 n \int_{0}^{1/a_K} \frac{dy}{(n + y)^2} \\
			&\leq b p_0 \| f'' \| \| \phi \|^2 \frac{1}{n a_K}.
		\end{align*}
		As a result, there exists a constant $ C > 0 $ such that, for all $ n \geq c_0/2 $,
		\begin{equation} \label{bound_RK2_neveu}
			|R^{f,\phi}_{K,2}(\rho, n)| \leq \frac{C}{a_K}.
		\end{equation}
		
		\paragraph*{Bound on $ R^{f,\phi}_{K,3} $}
		
		We start by writing, using \eqref{bound_Theta_small_y} in the first integral,
		\begin{multline*}
			|R^{f,\phi}_{K,3}(\rho, n)| \leq b \| f'' \| \| \phi \|^2 n \int_{1/a_K}^{n} \left( \frac{y_K}{n + y_K} \right)^2 \frac{| p(a_K y) - p_0 |}{y^2} dy \\ + 2 \| f \| b n \int_{n}^{\infty} \frac{| p(a_K y) - p_0 |}{y^2} dy.
		\end{multline*}
		Since $ y_K \leq y $, the first term on the right is bounded by
		\begin{equation*}
			b \| f'' \| \| \phi \|^2 \int_{1/(n a_K)}^{1} | p(n a_K y) - p_0 | dy. 
		\end{equation*}
		Changing variables in the second integral, we obtain that
		\begin{equation*}
			|R^{f,\phi}_{K,3}(\rho, n)| \leq h(n a_K),
		\end{equation*}
		where
		\begin{equation*}
			h(r) := b \| f'' \| \| \phi \|^2 \int_{1/r}^{1} | p(r y) - p_0 | dy + 2 \| f \| b \int_{1}^{\infty} \frac{| p(r y) - p_0 |}{y^2} dy.
		\end{equation*}
		By \eqref{condition_pk_neveu} and the dominated convergence theorem, $ h(r) \to 0 $ as $ r \to \infty $.
		As a result, there exists a sequence $ (\eta_K, K \geq 1) $ converging to zero as $ K \to \infty $ such that, for all $ n \geq c_0 / 2 $,
		\begin{equation} \label{bound_RK3_neveu}
			|R^{f,\phi}_{K,3}(\rho, n)| \leq \eta_K.
		\end{equation}
		
		\paragraph*{Bound on $ R^{f,\phi}_{K,4} $}
		
		To bound $ R^{f,\phi}_{K,4} $, we split the integral according to whether $ y \leq 1 $ or $ y > 1 $.
		When $ y > 1 $, we use \eqref{bound_f_yK_y}, yielding
		\begin{equation*}
			b p_0 n \int_{1}^{\infty} \abs{ \Theta_{n, y_K} F_{f,\phi}(\rho) - \Theta_{n, y} F_{f,\phi}(\rho) } \frac{dy}{y^2} \leq 2 b p_0 \| f' \| \| \phi \| \frac{1}{K} \int_{1}^{\infty} \frac{dy}{y^2}.
		\end{equation*}
		When $ y \leq 1 $, we use \eqref{bound_Theta_yK_y} to obtain
		\begin{align*}
			b p_0 n \int_{1/a_K}^{1} \abs{ \Theta_{n, y_K} F_{f,\phi}(\rho) - \Theta_{n, y} F_{f,\phi}(\rho) } \frac{dy}{y^2} &\leq \frac{2 C}{K} b p_0 \int_{1/a_K}^{1} \frac{dy}{y}, \\
			&= 2 C b p_0 \frac{\log(K \log(K))}{K}.
		\end{align*}
		As a result, there exists $ (\eta_K, K \geq 1) $ converging to zero as $ K \to \infty $ such that, for all $ n \geq c_0/2 $,
		\begin{equation} \label{bound_RK4_neveu}
			|R^{f,\phi}_{K,4}(\rho,n)| \leq \eta_K.
		\end{equation}
		
		\paragraph*{Conclusion of the proof}
		
		Combining \eqref{bound_RK1_neveu}, \eqref{bound_RK2_neveu}, \eqref{bound_RK3_neveu} and \eqref{bound_RK4_neveu}, we obtain that there exists a deterministic sequence $ (\eta_K, K \geq 1) $ converging to zero as $ K \to \infty $ such that, for all $ n \geq n_* / 2 $,
		\begin{equation*}
			\sum_{i=1}^{4} | R^{f,\phi}_{K,i}(\rho, n) | \leq \eta_K.
		\end{equation*}
		Then, setting
		\begin{equation*}
			\varepsilon_K^{f,\phi}(t) := -\int_{0}^{t \wedge \tau_K} \sum_{i=1}^{4} R^{f,\phi}_{K,i}(\rho^K_s, \Nbar(s)) ds,
		\end{equation*}
		we obtain the result, namely that 
		\begin{equation*}
			F_{f,\phi}(\rho^K_{t \wedge \tau_K}) - F_{f,\phi}(\rho^K_0) - \int_{0}^{t \wedge \tau_K} A F_{f,\phi}(\rho^K_s) ds + \varepsilon^{f,\phi}_K(t)
		\end{equation*}
		is an $ \F^K_t $-martingale and that
		\begin{equation*}
			\sup_{t \in [0, T]} \abs{ \varepsilon^{f,\phi}_K(t) } \leq \eta_K, \quad \text{ almost surely.}
		\end{equation*}
		This concludes the proof of the lemma.
	\end{proof}
	
	\subsection{Compact containment of the marginals of the scaled population size process} \label{subsec:proof_moment_bound_neveu}
	
	It remains to prove Lemma~\ref{lemma:bound_E_log_V_n}.
	To do so, we take $ L > C_0 $ and we define
	\begin{equation*}
		\tau_{K,L} := \inf \lbrace t \geq 0 : \Nbar(t) > L \rbrace \wedge \tau_K.
	\end{equation*}
	We then prove the following.
	
	\begin{lemma} \label{lemma:limsup_E_Nbar}
		For any $ L > n_* $ and any $ t \geq 0 $,
		\begin{equation*}
			\limsup_{K \to \infty} \E{ \int_{0}^{t \wedge \tau_{K,L}} \Nbar(s) ds } \leq t n_*.
		\end{equation*}
	\end{lemma}
	
	\begin{proof}
		Define a random measure $ \Gamma_K^{(L)} $ on $ \R_+ \times [c_0/2,L] $ by
		\begin{equation*}
			\Gamma_K^{(L)}([0,t] \times B) := \int_{0}^{t \wedge \tau_{K,L}} \1{\Nbar(s) \in B}  ds, \quad t \geq 0, B \in \mathcal{B}([c_0/2,L]).
		\end{equation*}
		By Lemma~1.1 in \cite{kurtz_averaging_1992}, $ (\Gamma_K^{(L)}, K \geq 1) $ is tight in $ (\ell([c_0/2,L]), d) $.
		Using Lemma~\ref{lemma:Nbar_neveu} and reasoning as in Example~2.3 in \cite{kurtz_averaging_1992}, we see that, for any limit point $ \Gamma^{(L)} $ of a converging subsequence of $ (\Gamma_K^{(L)}, K \geq 1) $ and for any $ \lambda > 0 $,
		\begin{equation*}
			\int_{[0,t] \times [c_0/2,L]} \mathcal{B} g_\lambda(y) \Gamma^{(L)}(ds, dy) = 0,
		\end{equation*}
		for all $ t \geq 0 $, almost surely.
		This implies that $ \Gamma^{(L)} $ is supported on $ \R_+ \times \lbrace n_* \rbrace $.

		To conclude the proof of the statement, we note that
		\begin{equation*}
			\int_{0}^{t \wedge \tau_{K,L}} \Nbar(s) ds = \int_{[0,t] \times [c_0/2,L]} y \, \Gamma_K^{(L)}(ds, dy).
		\end{equation*}
		Since $ (\Gamma_K^{(L)}, K \geq 1) $ is tight, from any subsequence we can extract a converging subsequence.
		Along such a subsequence,
		\begin{equation*}
			\int_{[0,t] \times [c_0/2,L]} y \, \Gamma_K^{(L)}(ds, dy) \to \int_{[0,t] \times [c_0/2,L]} y \, \Gamma^{(L)}(ds, dy),
		\end{equation*}
		in distribution.
		Moreover,
		\begin{equation*}
			\int_{[0,t] \times [c_0/2,L]} y \, \Gamma_K^{(L)}(ds, dy) \leq L t,
		\end{equation*}
		almost surely.
		Thus, by the dominated convergence theorem, along such a converging subsequence,
		\begin{equation*}
			\E{\int_{[0,t] \times [c_0/2,L]} y \, \Gamma_K^{(L)}(ds, dy)} \to \E{\int_{[0,t] \times [c_0/2,L]} y \, \Gamma^{(L)}(ds, dy)}.
		\end{equation*}
		Since $ \Gamma^{(L)} $ is supported on $ \R_+ \times \lbrace n_* \rbrace $, the right hand side is smaller than $ n_* t $.
		As a result, for any $ \varepsilon > 0 $, one cannot find a subsequence along which
		\begin{equation*}
			\E{\int_{0}^{t \wedge \tau_{K,L}} \Nbar(s) ds} \geq n_* t + \varepsilon,
		\end{equation*}
		which concludes the proof of Lemma~\ref{lemma:limsup_E_Nbar}.
	\end{proof}
	
	Recall the definition of $ \overline{n}_K $ in \eqref{def:n_bar_K}.
	We define $ V_K : \R_+ \to \R_+ $ as
	\begin{equation*}
		V_K(n) := \frac{n}{\overline{n}_K} - 1 - \log\left( \frac{n}{\overline{n}_K} \right).
	\end{equation*}
	We also set
	\begin{equation*}
		g_K(n) := \log(1 + V_K(n)),
	\end{equation*}
	and, for any $ \lambda > 0 $,
	\begin{equation*}
		f_{\lambda,K}(n) := \frac{\lambda g_K(n)}{\lambda + g_K(n)}.
	\end{equation*}
	By \eqref{cvg_n_bar_K}, for any $ \varepsilon \in (0, n_*) $, $ V_K \geq V_\varepsilon $ for all $ K $ large enough, so Lemma~\ref{lemma:bound_E_log_V_n} will be proved if we show that there exists $ C_{T} > 0 $ such that
	\begin{equation} \label{bound_E_gK}
		\sup_{K \geq 1} \sup_{t \in [0,T]} \E{ g_K(\Nbar(t \wedge \tau_K)) } \leq C_{T}.
	\end{equation}
	We shall instead work with $ f_{\lambda, K} $ as it is bounded, letting $ \lambda \to \infty $ at the end to obtain the above bound.
	
	The first step will be to prove the following.
	Recall the definition of the operator $ \mathcal{B}_K $ in \eqref{def:BK_neveu}.
	
	\begin{lemma} \label{lemma:bound_B_f_lambda_K}
		There exist constants $ \newCst{fl1} > 0 $ and $ \newCst{fl2} > 0 $ such that, for all $ \lambda > 0 $ and for any $ n \geq \frac{c_0}{2} $,
		\begin{equation*}
			\mathcal{B}_K f_{\lambda,K} (n) \leq \Cst{fl1} n + \Cst{fl2}.
		\end{equation*}
	\end{lemma}
	
	\begin{proof}
		Using \eqref{decomp:BK_neveu}, we write
		\begin{equation} \label{decomposition_B_f}
			\mathcal{B}_K f_{\lambda, K} (n) = c \log(K) (\overline{n}_K - n) n f_{\lambda, K}'(n) + \sum_{i=2}^{4} B^K_i f_{\lambda,K}(n).
		\end{equation}
		Since
		\begin{equation*}
			f_{\lambda,K}'(n) = \left( \frac{\lambda}{\lambda + g_K(n)} \right)^2 \frac{1}{1 + V_K(n)} \left( \frac{1}{\overline{n}_K} - \frac{1}{n} \right),
		\end{equation*}
		we see that
		\begin{equation*}
			(\overline{n}_K - n) n f_{\lambda, K}'(n) \leq 0.
		\end{equation*}
		We then bound each $ B^K_i f_{\lambda,K}(n) $ separately.
		
		\paragraph*{Bound on $ B^K_2 f_{\lambda,K}(n) $}
		
		To bound $ B^K_2 f_{\lambda,K}(n) $, we note that, for all $ y \geq 0 $ and all $ n \geq 0 $,
		\begin{align*}
			V_K(n + y) - V_K(n) \leq \frac{y}{\overline{n}_K} \leq \frac{y}{n_* - \varepsilon}
		\end{align*}
		for all $ K $ large enough, by \eqref{cvg_n_bar_K}.
		In addition, for any $ \delta > 0 $ and $ V > 0 $,
		\begin{align*}
			\log(1+V+\delta) - \log(1+V) &= \log\left( 1 + \frac{\delta}{1 + V} \right) \\
			&\leq \log(1+\delta).
		\end{align*}
		As a result, for $ n \geq 0 $ and $ y \geq 0 $,
		\begin{equation*}
			f_{\lambda,K}(n + y) - f_{\lambda,K}(n) \leq \log\left( 1 + \frac{y}{n_* - \varepsilon} \right).
		\end{equation*}
		Plugging this in the definition of $ B^K_2 f_{\lambda,K} $ in \eqref{decomp:BK_neveu}, we obtain
		\begin{equation*}
			B^K_2 f_{\lambda,K}(n) \leq a_K b n \sum_{k > a_K} p_k \log\left( 1 + \frac{k}{a_K (n_* - \varepsilon)} \right).
		\end{equation*}
		Using \eqref{def:p_y}, we then write
		\begin{align*}
			B^K_2 f_{\lambda,K}(n) &\leq a_K b n \int_{\lfloor a_K \rfloor + 1}^{\infty} \log\left( 1 + \frac{\lfloor y \rfloor}{a_K (n_* - \varepsilon)} \right) \frac{p(y)}{y^2} dy \\
			&= b n \int_{\frac{\lfloor a_K \rfloor + 1}{a_K}}^{\infty} \log\left( 1 + \frac{\lfloor a_K y \rfloor / a_K}{n_* - \varepsilon} \right) \frac{p(a_K y)}{y^2} dy \\
			&\leq b \, C_p \, n \int_{1}^{\infty} \log\left( 1 + \frac{y}{n_* - \varepsilon} \right) \frac{dy}{y^2},
		\end{align*}
		where $ C_p = \sup_{y \in [1,\infty)} p(y) $ was defined in \eqref{def:Cp}.
		Since the integral on the right hand side is finite, we obtain that there exists a constant $ C > 0 $ such that
		\begin{equation} \label{bound_B2_flambda}
			B^K_2 f_{\lambda,K}(n) \leq C n.
		\end{equation}
		
		\paragraph*{Bound on $ B^K_3 f_{\lambda,K}(n) $}
		
		By the definition of $ B^K_3 $,
		\begin{equation} \label{Taylor_B3}
			B^K_3 f_{\lambda,K}(n) = \frac{1}{a_K} b n \sum_{k=1}^{\lfloor a_K \rfloor} k^2 p_k \int_{0}^{1} (1-u) f_{\lambda,K}''\left( n + u\frac{k}{a_K} \right) du.
		\end{equation}
		We then note that
		\begin{equation*}
			f_{\lambda,K}''(n) = \left( \frac{\lambda}{\lambda + g_K(n)} \right)^2 g_K''(n) - 2 \frac{\lambda^2}{(\lambda + g_K(n))^3} (g_K'(n))^2,
		\end{equation*}
		and
		\begin{equation*}
			g_K''(n) = - \frac{1}{(1+V_K(n))^2} \left( \frac{1}{\overline{n}_K} - \frac{1}{n} \right)^2 + \frac{1}{1 + V_K(n)} \times \frac{1}{n^2}.
		\end{equation*}
		As a result,
		\begin{align*}
			f_{\lambda,K}''(n) &\leq \left( \frac{\lambda}{\lambda + g_K(n)} \right)^2 \frac{1}{1 + V_K(n)} \times \frac{1}{n^2} \\
			&\leq \frac{1}{n^2}. \numberthis \label{bound_flambda_seconde}
		\end{align*}
		Plugging this in \eqref{Taylor_B3} and using the fact that $ \sup_{k \in \N} k^2 p_k < \infty $ by \eqref{condition_pk_neveu}, we obtain that there exists a constant $ C > 0 $ such that, for all $ n \geq \frac{c_0}{2} $,
		\begin{equation} \label{bound_B3_flambda}
			B^K_3 f_{\lambda,K}(n) \leq C
		\end{equation}
		
		\paragraph*{Bound on $ B^K_4 f_{\lambda,K}(n) $}
		
		By the definition of $ B^K_4 $,
		\begin{equation*}
			B^K_4 f_{\lambda,K}(n) = \frac{1}{a_K} (d + c \log(K) n) n \int_{0}^{1} (1-u) f_{\lambda,K}''\left( n - \frac{u}{a_K} \right) du.
		\end{equation*}
		Using \eqref{bound_flambda_seconde} again, we obtain that
		\begin{equation*}
			B^K_4 f_{\lambda,K}(n) \leq \frac{1}{K} \left( \frac{d}{\log(K)} + c n \right) \frac{n}{(n - 1/a_K)^2}.
		\end{equation*}
		As a result, there exists a constant $ C > 0 $ such that, for all $ n \geq \frac{c_0}{2} $ and for all $ K $ large enough,
		\begin{equation} \label{bound_B4_flambda}
			B^K_4 f_{\lambda, K}(n) \leq \frac{C}{K}.
		\end{equation}
		
		\paragraph*{Conclusion of the proof}
		
		Plugging \eqref{bound_B2_flambda}, \eqref{bound_B3_flambda} and \eqref{bound_B4_flambda} in \eqref{decomposition_B_f}, we obtain that, for all $ n \geq \frac{c_0}{2} $,
		\begin{equation*}
			\mathcal{B}_K f_{\lambda, K}(n) \leq \Cst{fl1} n + \Cst{fl2},
		\end{equation*}
		which concludes the proof of the lemma.
	\end{proof}
	
	We are now able to prove Lemma~\ref{lemma:bound_E_log_V_n}.
	
	\begin{proof}[Proof of Lemma~\ref{lemma:bound_E_log_V_n}]
		Since $ n \mapsto f_{\lambda,K}(n) $ is bounded, for any $ L > C_0 $, $ \lambda > 0 $ and $ t \geq 0 $, by Lemma~\ref{lemma:bound_B_f_lambda_K},
		\begin{equation*}
			\E{ f_{\lambda,K}(\Nbar(t \wedge \tau_{K,L})) } \leq \E{ f_{\lambda,K}(\Nbar(0)) } + \Cst{fl1} \E{ \int_{0}^{t \wedge \tau_{K,L}} \Nbar(s) ds } + \Cst{fl2} t.
		\end{equation*}
		Since $ f_{\lambda,K}(n) \leq g_K(n) $ for all $ n \geq 0 $, for all $ \lambda > 0 $,
		\begin{equation*}
			\E{ f_{\lambda,K}(\Nbar(0)) } \leq \E{ g_K(\Nbar(0)) }.
		\end{equation*}
		By \eqref{N0_compact_neveu} and \eqref{cvg_n_bar_K}, the right hand side is bounded uniformly in $ K $.
		By Lemma~\ref{lemma:limsup_E_Nbar}, we obtain that, for any $ T \geq 0 $, there exists a constant $ C_T > 0 $ such that, for all $ t \in [0,T] $, $ \lambda > 0 $, $ L > n_* $,
		\begin{equation*}
			\E{ f_{\lambda,K}(\Nbar(t\wedge \tau_{K,L})) } \leq C_T.
		\end{equation*}
		Using the fact that $ n \mapsto f_{\lambda,K}(n) $ is bounded, by the dominated convergence theorem, we can take the limit as $ L \to \infty $ to obtain
		\begin{equation*}
			\sup_{t \in [0,T]} \E{ f_{\lambda,K}(\Nbar(t \wedge \tau_K)) } \leq C_T.
		\end{equation*}
		By the monotone convergence theorem, we can take the limit as $ \lambda \to \infty $ to obtain
		\begin{equation*}
			\sup_{K \geq 1} \sup_{t \in [0,T]} \E{ g_K(\Nbar(t\wedge \tau_K)) } \leq C_T.
		\end{equation*}
		By \eqref{bound_E_gK}, this concludes the proof of Lemma~\ref{lemma:bound_E_log_V_n}.
	\end{proof}

	\bibliography{birth-death}
\end{document}